\documentclass[11pt,a4paper]{deepmind}
\pdftrailerid{redacted}

\makeatletter
\renewcommand\bibentry[1]{\nocite{#1}{\frenchspacing\@nameuse{BR@r@#1\@extra@b@citeb}}}
\makeatother

\usepackage{kantlipsum, lipsum}
\usepackage{dsfont}
\usepackage{dm-colors}
\usepackage{amsthm} 
\usepackage[normalem]{ulem}
\usepackage{multirow}
\usepackage{algorithm}
\usepackage{algorithmic}
\usepackage{mathtools}
\usepackage{subfigure}
\usepackage{subcaption}
\usepackage{graphicx}

\usepackage{cleveref}
\usepackage[toc,page]{appendix}

\usepackage{threeparttable}

\makeatletter
  \def\vhrulefill#1{\leavevmode\leaders\hrule\@height#1\hfill \kern\z@}
\makeatother

\usepackage{tikz} 
\usepackage{eso-pic}

\newcommand\BackgroundPicture{%
  \put(0,0){%
    \parbox[b][\paperheight]{\paperwidth}{%
      \vfill
      \centering%
\begin{tikzpicture}[remember picture,overlay]
\node [rotate=60,scale=8,text opacity=0.2] at (current page.center) {}; 
\end{tikzpicture}%
      \vfill
    }}}

\usepackage[numbers, sort&compress, square, round]{natbib}

\graphicspath{{figures/}}

\title{
Exploiting Electrolyzer Flexibility via Multiscale Model Predictive Control Cross Heterogeneous Energy Markets
}

\correspondingauthor{Jiaze Ma}

\reportnumber{} 

\author[1]{Zhichao Chen}
\author[2]{Hongyuan Sheng}
\author[3]{Hao Wang}
\author[1]{Jiaze Ma}

\affil[1]{Department of Systems Engineering, City University of Hong Kong\\}
\affil[2]{Department of Chemistry, Fudan University\\}
\affil[3]{Department of Automation, Zhejiang University}

\begin{abstract}
Green hydrogen production via electrolysis is crucial for decarbonization but faces significant economic hurdles primarily due to the high cost of the electricity. 
However, current electrolyzer-based hydrogen production processes predominantly rely on the single-scale Day-Ahead Market (DAM) for electricity procurement, failing to fully exploit the economic benefits offered by multi-scale electricity market that integrates both the DAM and the Real-Time Market (RTM), thereby eliminating the opportunity to reduce the overall cost.
To mitgate this technical gap, this research investigates a dynamic operational strategy enabling electrolyzers to strategically navigate between the DAM and RTM to minimize net operation costs.
Using an rolling horizon optimization framework to coordinate bidding and operation, we demonstrate a strategy where electrolyzers secure primary energy via exclusive DAM purchases, then actively engage the RTM to buy supplemental energy cheaply or, critically, sell procured DAM energy back at a profit during high RTM price periods.
Our analysis reveals that this coordinated multi-scale electricity market participation strategy can dramatically reduce net electricity expenditures, achieving near-zero or even negative effective electricity costs for green hydrogen production under realistic market scenarios, effectively meaning the operation can profit from its electricity market interactions.
By transforming electrolyzers from simple price-takers into active participants capable of arbitrage between market timescales, this approach unlocks a financially compelling pathway for green hydrogen, accelerating its deployment while simultaneously enhancing power grid flexibility. 

\end{abstract}

\begin{document}
\AddToShipoutPicture{\BackgroundPicture}
\newtheorem{theorem}{Theorem}[section] 
\newtheorem{definition}[theorem]{Definition} 
\newtheorem{lemma}{Lemma} 
\newtheorem{corollary}[theorem]{Corollary}
\newtheorem{example}{Example}[section]
\newtheorem{proposition}[theorem]{Proposition}

\maketitle


\section{Introduction}\label{sec:intro}
The global imperative to transition to a net-zero carbon emission system is under immense pressure, as emphasized by the binding international climate agreements and national mandates demanding rapid, cross-sectoral decarbonization~\cite{erickson2023biogas,sinha2024diverse}. Green hydrogen has been identified as a critical enabler in this transition, especially for hard-to-abate industries where direct electrification is challenging~\cite{yang2022breaking,guan2023hydrogen,wallington2024green,shafiee2024carbon,segovia2025green}. However, large-scale adoption fundamentally depends on achieving cost competitiveness with incumbent fossil-fuel-based hydrogen and alternative low-carbon options such as blue hydrogen~\cite{alhumaidan2023blue,ueckerdt2024cost}. This economic imperative is further amplified by concerns over energy security and the strategic desire to leverage abundant domestic renewable resources~\cite{nishiyama2021photocatalytic,sakib2024harnessing}. Furthermore, the anticipation of future carbon pricing mechanisms is expected to incentivize proactive shifts toward affordable clean energy~\cite{zantye2021renewable,jiang2024unclean}. Thus, developing innovative operational approaches and market interaction strategies to drastically reduce the overall cost for green hydrogen manufacture remains a paramount challenge for the global energy transition.

The core technology enabling green hydrogen is water electrolysis, a process that utilizes electricity to split water ($\text{H}_2\text{O}$) into hydrogen ($\text{H}_2$) and oxygen ($\text{O}_2$)~\cite{ursua2011hydrogen,chatenet2022water,shi2023customized,lee2023efficient,shi2023sodium,wang2024water}, and among the various electrolyzer technologies, Proton Exchange Membrane (PEM) electrolyzers are particularly well-suited for dynamic operation due to their rapid response capabilities and high operational flexibility. However, while this operational flexibility does create the potential for significant economic savings, it simultaneously introduces a critical trade-off. On one hand, the overall cost of hydrogen production is heavily influenced by the electricity expenditure, as electricity accounts for the largest share of the operational cost. On the other hand, frequent and dynamic operation—though beneficial for leveraging favorable electricity prices—can substantially accelerate the degradation of the PEM membrane, which is one of the most critical and costly components of the PEM-based hydrogen production process. Thus, the economic viability of PEM-based green hydrogen production ultimately hinges on carefully \emph{balancing} electricity cost minimization with prudent management of membrane replacement costs.

\begin{figure}[!h]
    \centering
    \includegraphics[width=1\textwidth]{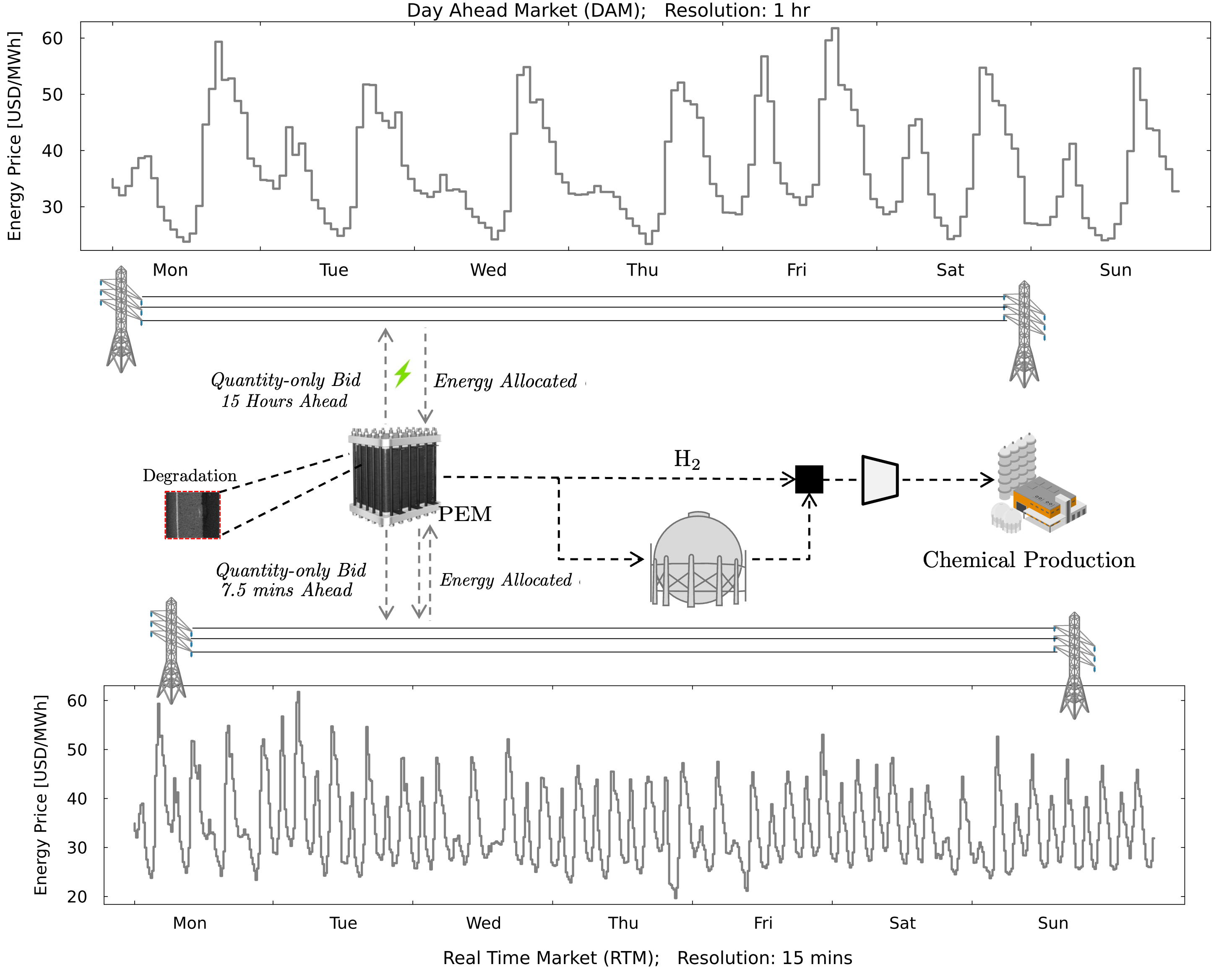}   
    \caption{Bidding from the DAM and RTM with a flexible PEM.}\label{fig:electrolyzerProcess}
\end{figure}

Fortunately, the accelerating integration of variable renewable energy sources, such as wind and solar, is fundamentally reshaping electricity markets, particularly by increasing price volatility and creating significant spreads between Day-Ahead Market (DAM) and Real-Time Market (RTM)~\cite{8244305,gao2022multiscale,cao2020multiscale} (see \Cref{fig:electrolyzerProcess}). In light of these evolving circumstances, modern market frameworks are increasingly enabling large and flexible consumers—such as electrolyzer-based hydrogen facilities—to move beyond the passive role of price-takers and assume a more active, participatory position in electricity markets. By leveraging temporal price fluctuations through optimized bidding and dynamic adjustment of electricity consumption in both DAM and RTM, these consumers are now positioned to exploit arbitrage opportunities and unlock additional value streams, such as demand response and grid balancing services~\cite{cao2020multiscale}. Such capabilities stand in marked contrast to traditional procurement approaches, like fixed-price power purchase agreements, which by design cannot capture the value of operational flexibility under volatile market conditions~\cite{falcao2020review,schofield2025dynamic}.

The inherent operational flexibility of PEM electrolyzers—characterized by rapid ramping and wide operating ranges—makes them especially suitable for market-driven participation. However, translating this technical potential into practical economic benefits is not straightforward. While flexible operation enables cost savings and revenue enhancement, it also introduces complex technical challenges: Specifically, frequent power cycling and variable load profiles associated with active market participation can accelerate membrane degradation. The economic viability of PEM-based hydrogen production therefore critically depends on not only maximizing market revenues, but also minimizing lifecycle costs associated with equipment degradation and replacement~\cite{feng2017review,tomic2023review,norazahar2024degradation}. Addressing these challenges requires the development of sophisticated optimization frameworks that holistically balance electricity procurement costs against the degradation-induced expenses of membrane replacement. In this context, integrating high-fidelity degradation models into the operational strategy for PEM electrolyzers becomes essential. Such a modeling approach makes it possible to rigorously evaluate the trade-offs inherent to flexible market participation and devise operational schedules that achieve both cost competitiveness and sustainable asset longevity.

To alleviate the abovementioned issues, this manuscript presents and rigorously evaluates an optimized operational strategy that enables green hydrogen facilities to actively and intelligently participate across both the DAM and RTM (see~\Cref{fig:electrolyzerProcess} for a model illustration). The proposed approach is based on the rolling horizon optimization framework, which is fully integrated with a high-fidelity PEM electrolyzer model~\cite{chandesris2015membrane} that accurately captures the complex system behaviros of membrane degradation under fluctuating operational profiles. The rolling horizon optimization co-optimizes both market bidding quantities and operational schedules across multiple market scales, explicitly balancing the pursuit of dynamic electricity market revenues and cost savings against the projected costs of membrane degradation and replacement~\cite{schofield2025dynamic,cao2020multiscale}. Throughout this study, we adopt the price taker assumption similar to previous work~\cite{lu2004pumped,kanakasabapathy2010bidding,krishnamurthy2017energy} for the hydrogen facilities, in which market participants optimize their bids and schedules based on the given prices without exerting market power. For both DAM and RTM, we assume that the hydrogen facility submits only quantity bids (i.e., no offer price is specified), and is able to transact any requested volume at the prevailing market price. In the DAM, only electricity purchases are allowed (i.e., quantity to be bought is bid for 24 hours in advance). In the RTM, the facility can adjust its position by both buying and selling electricity. Financial settlements are performed at the market clearing prices for DAM and RTM, respectively. The efficacy of the proposed approach is demonstrated through validation with historical electricity price data from the Houston and California markets. 


Our experimental results reveal that integrating high-fidelity membrane degradation model with multi-scale market participation can dramatically improve the economic performance of green hydrogen production. Specifically, this approach consistently results in substantial reductions in the levelized cost of hydrogen, and in certain market conditions, can even yield negative incremental operational cost periods for PEM electrolyzers—reflecting net profits from strategic electricity market interactions. More importantly, our findings expose several counter-intuitive operational advantages that run counter to conventional industry intuitions and reveal key insights into the economic and technical optimization of PEM-based green hydrogen production. In particular, we find that:
\begin{itemize}
    \item{Contrary to the expectation that aggressive operational changes might accelerate membrane degradation, employing the high-fidelity operational model even within a single scale electricity market (HF-SS) not only substantially curtails electricity expenditure compared to traditional constant operation, but also, significantly mitigates membrane degradation. This demonstrates that optimized, dynamic strategy, guided by accurate system models, can concurrently achieve energy cost savings and extend the lifetime of capital-intensive membrane.}

\item{Extending this high-fidelity dynamic approach to strategically navigate multiple electricity market scales (HF-MS)—specifically co-optimizing across DAM and RTM—unlocks further substantial reductions in the overall levelized cost of hydrogen. This superiority over even optimized HF-SS underscores the significant economic value captured by leveraging inter-market price differentials and volatilities through arbitrage.}

\item{Within the complex dynamics of multi-scale market participation, the precision of the underlying system model is paramount. Our results show that the HF-MS achieves demonstrably superior outcomes in terms of minimizing membrane degradation and reducing total cumulative operational costs when compared to the multi-scale strategy reliant on the lower-fidelity model (LF-MS). This highlights that accurately capturing and predicting electrolyzer behavior, especially membrane degradation, is crucial to fully exploit the economic potential of multi-market operations without incurring prohibitive membrane degradation costs.}
\end{itemize}

\section{Experimental Results}

This section presents the comparative performance of the proposed HF-MS trading strategy against benchmark strategies within the Houston electricity market (the results for California electricity market are given in the Appendix). We evaluate its efficacy in market arbitrage, its impact on operational patterns affecting membrane degradation, and the resultant overall economic benefits for hydrogen production.

\begin{figure}[!h]
 \centering
    \includegraphics[width=0.95\textwidth]{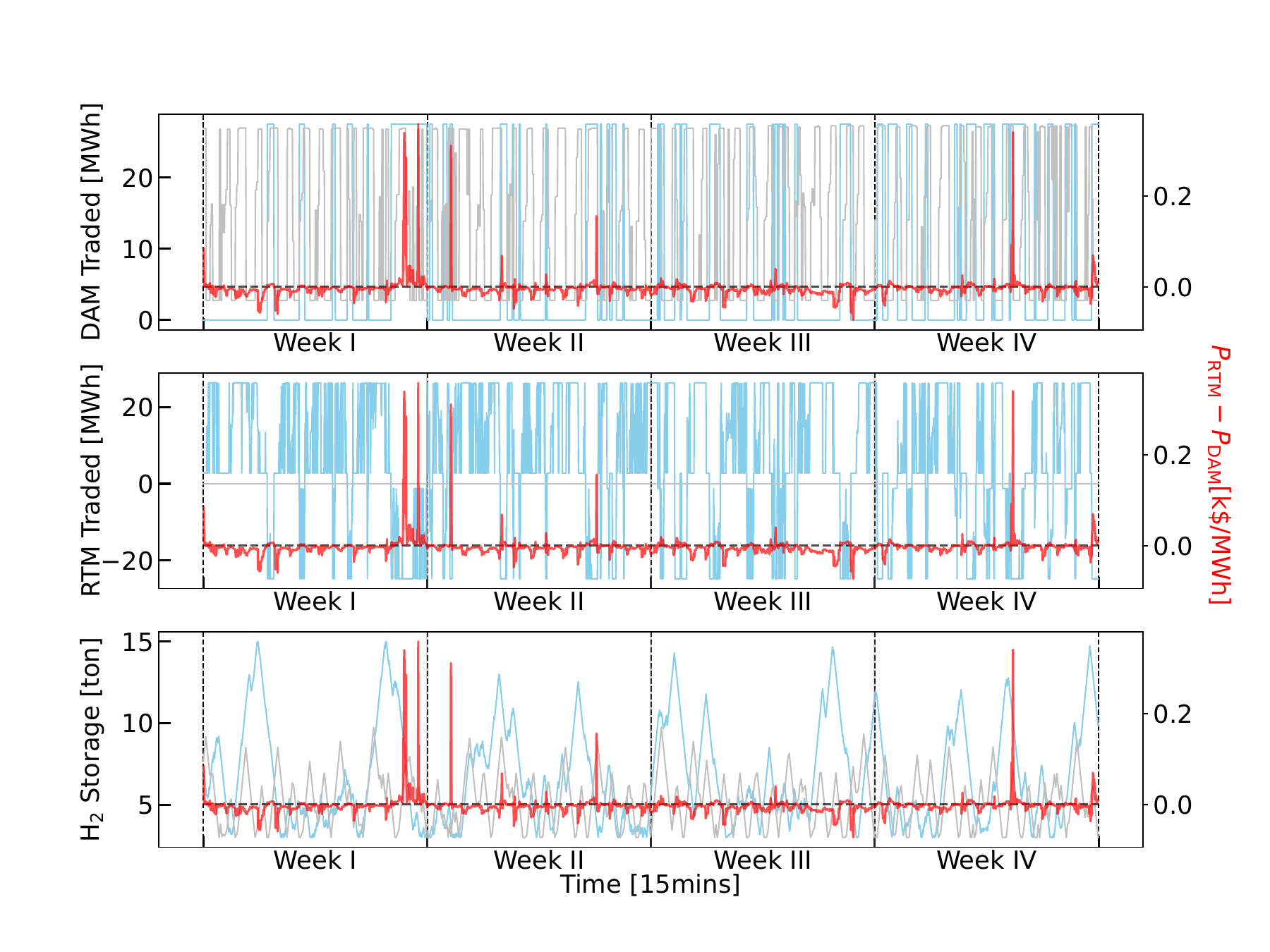}
    \caption{DAM traded power, RTM traded power, and $\text{H}_{\text{2}}$ storage levels for the HF-MS (light blue line) and HF-SS (gray line) strategies. Data correspond to four weeks in Houston, January 2022. The red line indicates the price differential (RTM Price minus DAM Price).}
    \label{fig:market_simple_ms_ss}
\end{figure}

\Cref{fig:market_simple_ms_ss} illustrates the distinct operational strategy of the HF-MS approach compared to the HF-SS strategy over four representative weeks. The operational data for HF-MS and HF-SS are detailed in the figure, where the top subplot displays traded DAM power (blue for HF-MS, gray for HF-SS), the middle subplot shows the traded RTM power (blue for HF-MS, gray for HF-SS), and the bottom subplot depicts hydrogen storage levels (blue for HF-MS, gray for HF-SS), and the red line indicates the electricity price differences between RTM and DAM.

Several key operational patterns for HF-MS emerge. Firstly, regarding DAM interaction (top subplot of \Cref{fig:market_simple_ms_ss}(a)), while HF-MS may engage in DAM trading with a different frequency profile compared to HF-SS, it strategically commits to significant DAM power purchases, particularly when anticipating or observing high RTM prices (indicated by spikes in the red line). This proactive DAM procurement aims to secure lower-cost electricity. Secondly, building upon these DAM purchases, the HF-MS strategy demonstrates an ability to capitalize on favorable RTM conditions. As seen in of \Cref{fig:market_simple_ms_ss}(b), during periods of high RTM prices (red line), HF-MS is observed to sell power back to the RTM (indicated by negative values for RTM Power), thereby realizing profits from the price arbitrage between DAM purchases and RTM sales. For example, prominent instances of selling power into the RTM during price spikes can be observed in Week I and Week IV. Finally, the operational flexibility of the HF-MS strategy is evident in the hydrogen storage dynamics (bottom subplot of \Cref{fig:market_simple_ms_ss}(c)). The $\text{H}_{\text{2}}$ storage levels exhibit frequent and considerable fluctuations. This behavior underscores the HF-MS's capacity for dynamic adjustment of hydrogen production in response to real-time market conditions and energy arbitrage opportunities, leading to more intensive utilization and regulation of the electrolysis plant.

\begin{figure}[!h]
 \centering
    \includegraphics[width=0.95\textwidth]{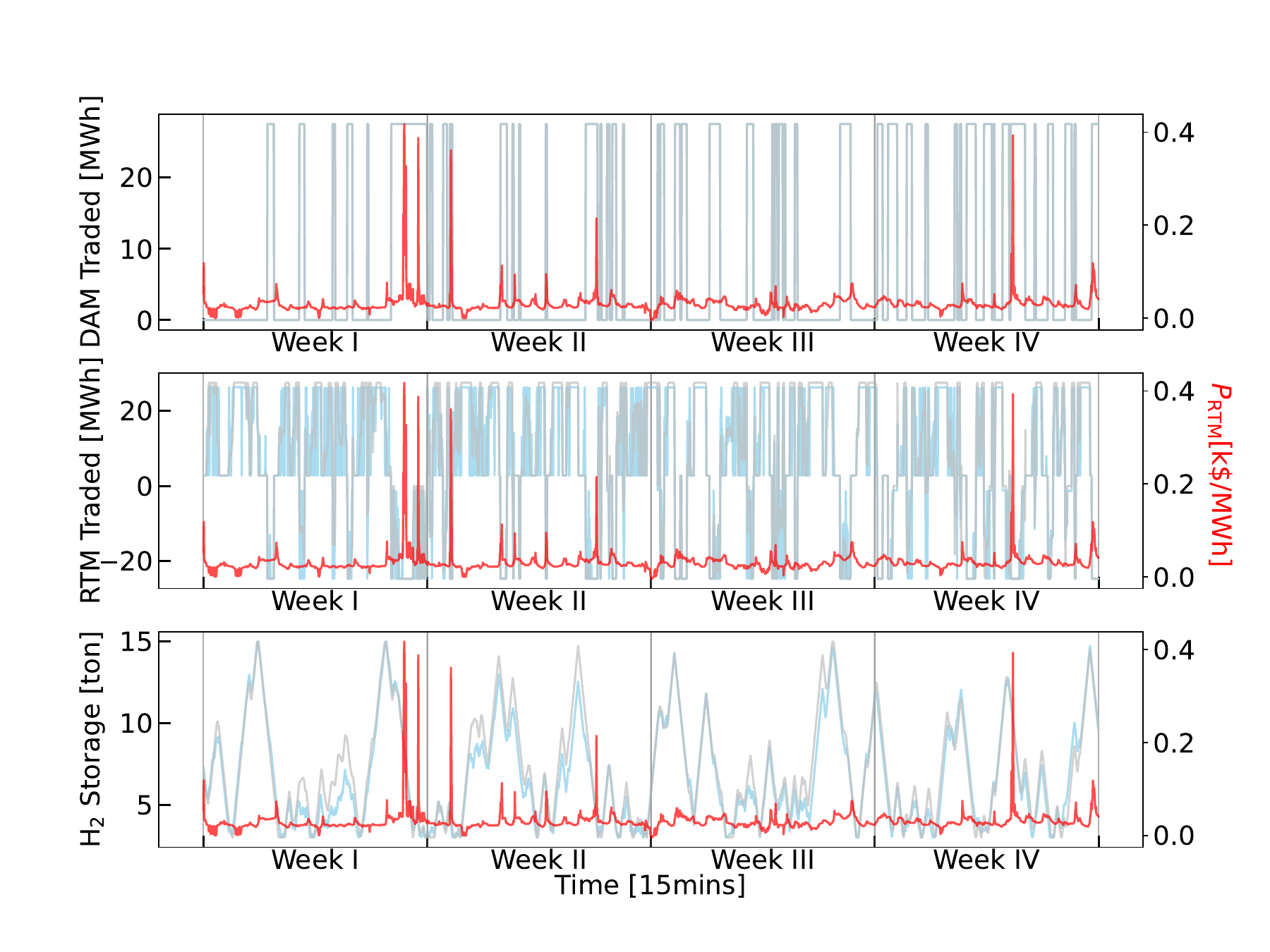}
        \caption{DAM traded power, RTM traded power, and $\text{H}_{\text{2}}$ storage levels for the HF-MS (light blue line) and LF-MS (gray line) strategies. Data correspond to four weeks in Houston, January 2022. The red line indicates the price differential (RTM Price minus DAM Price).}
    \label{fig:market_simple_hf_lf}
\end{figure}
Further comparative analysis was conducted between the HF-MS and LF-MS strategies, with HF-MS's key operational metrics depicted by the light blue lines in \Cref{fig:market_simple_hf_lf}. While both strategies may employ analogous fundamental logic for engaging with DAM and RTM—primarily exploiting price differentials for arbitrage—their specific interaction patterns within the RTM, particularly the magnitude of transactions and resultant systemic impacts, reveal notable distinctions.

A salient observation is the comparatively intensified RTM engagement exhibited by the LF-MS strategy. The middle subplot of \Cref{fig:market_simple_hf_lf}, which details RTM power transactions, clearly illustrates that when LF-MS (gray line) purchase power from the RTM (indicated by negative RTM Power values), the peak magnitudes and volumes of these sales are often more substantial. For instance, during the RTM price drop in Week I to Week IV (red line indicating RTM price), the LF-MS strategy not only participates in purchasing power but the extent of these sales (depth of the light gray line beyond zero) appears particularly pronounced, suggesting a propensity for larger-volume reverse power flows when profitable. Consequentially, this more aggressive RTM interaction under the LF-MS strategy correlates with markedly more pronounced fluctuations in hydrogen storage levels (bottom subplot, gray line). The $\text{H}_{\text{2}}$ storage undergoes cycles of more substantial accumulation. This behavior indicates that LF-MS's more pronounced power purchasing activity in the RTM—often undertaken to capitalize on favorable pricing for intensified hydrogen production—directly translates to more frequent and larger-amplitude variations in stored hydrogen. Such intensified procurement efforts, leading to accelerated hydrogen generation and storage cycles, reflect a more dynamic, and potentially more strenuous, utilization of the hydrogen storage capacity.

\begin{figure}[!h]
 \centering
    \subfigure[Cumulative cost of Houston under various operation strategies.]{
    \includegraphics[width=1.00\textwidth]{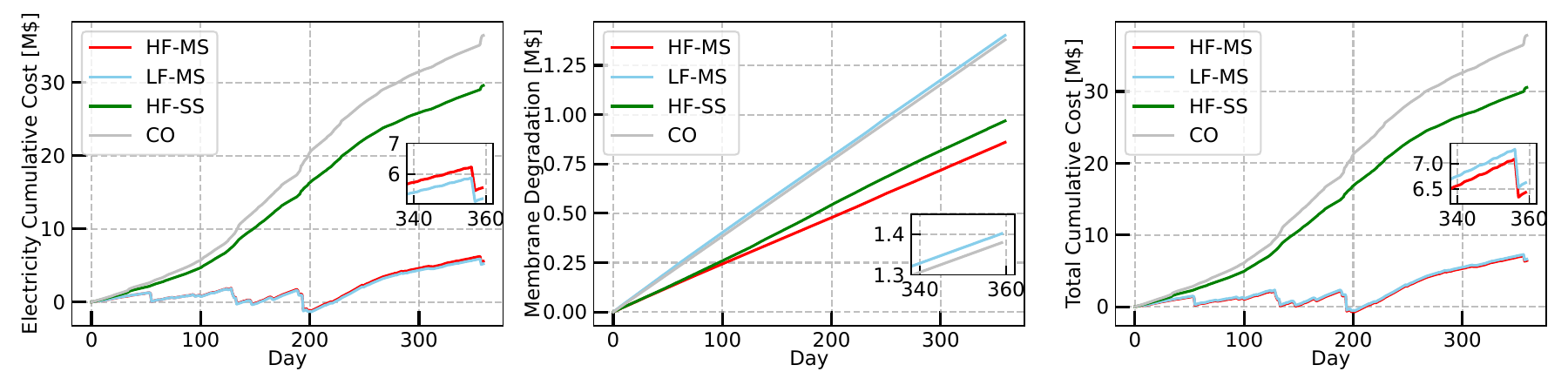}\label{subfig:cumCostVariousOperation}
} 
    \subfigure[Cumulative cost attribution for different operation strategies of Houston.]{
    \includegraphics[width=1.00\textwidth]{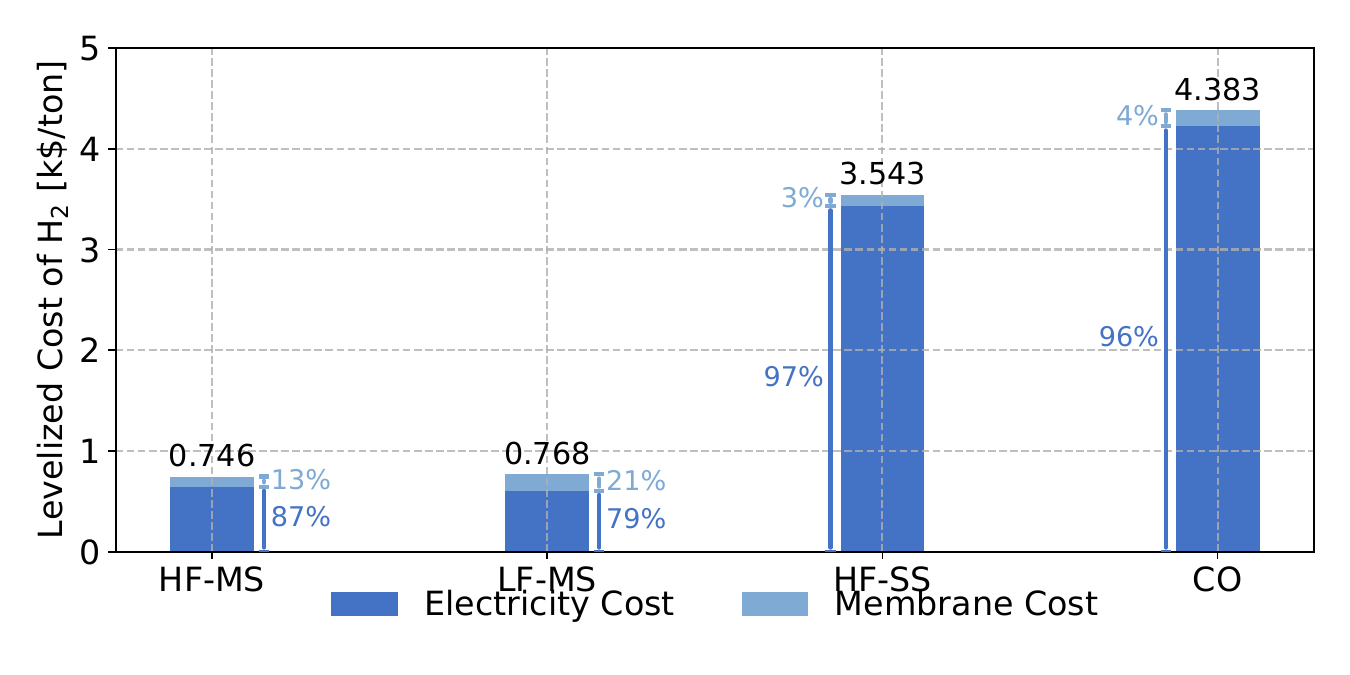}\label{subfig:cumCostbreakUp}
} 
    \caption{Cost comparison results.}
    \label{fig:market_cost_compare}
\end{figure}
To further clarify the economic impact of the proposed operational strategies, we compare the cumulative costs associated with various operationa stragies—namely, HF-MS, LF-MS, HF-SS, and constant operation (CO)—as illustrated in~\Cref{subfig:cumCostVariousOperation}. Building on this analysis, we decompose the resulting levelized cost of hydrogen (in k\$/ton) into its primary components—electricity consumption and membrane degradation—over the entire year for each strategy, as shown in~\Cref{subfig:cumCostbreakUp}.

From~\Cref{subfig:cumCostVariousOperation} the following insights can be given:
\begin{enumerate}
    \item{The traditional CO strategy, which runs the PEM electrolyzer at a fixed output based on DAM prices without any dynamic adjustment, leads to the highest cumulative costs—both in terms of electricity consumption and membrane degradation. This highlights the substantial economic penalty of not utilizing any form of optimization or operational flexibility.}
\item{Introducing a high-fidelity PEM model into the operation decision process within a single market (HF-SS) greatly reduces both electricity and membrane degradation costs compared to CO. Notably, jointly optimizing for cost and degradation not only lowers operational expenses, but also extends membrane life—demonstrating the essential value of high-fidelity model-based operation decision over naive constant operation. }
\item{Allowing market participation flexibility across DAM and RTM with LF-MS significantly reduces cumulative electricity costs compared to single-scale operation. However, the absence of precise membrane degradation modeling in LF-MS leads to much higher cumulative degradation costs, offsetting some of the gains in electricity savings. In contrast, HF-MS—combining market flexibility and high-fidelity degradation modeling—consistently delivers the lowest total cumulative cost. This underscores that capturing and accurately managing degradation is critical for fully realizing the value of multi-scale electricity market arbitrage.}
\end{enumerate}

Based on~\Cref{fig:market_cost_compare}, several notable insights can be drawn as follows:
\begin{enumerate}
    \item{Both HF-MS and LF-MS strategies achieve remarkably low levelized hydrogen production costs (0.746 and 0.768 k\$/ton, respectively), in stark contrast to the much higher costs observed in HF-SS (3.543 k\$/ton) and CO (4.383 k\$/ton) scenarios.
    }
\item{
In the HF-MS and LF-MS cases, electricity expenditures account for the majority of the levelized cost of hydrogen (87\% and 79\%, respectively), while membrane degradation remains a minor component (13\% and 21\%). Conversely, under the HF-SS and CO strategies, membrane degradation dominates the total cost, comprising 97\% and 96\%, respectively. This stark contrast highlights the critical impact of membrane wear when operation is restricted to DAM participation alone, underscoring the importance of both model fidelity and market flexibility in minimizing overall costs.
}
\item{
The HF-MS strategy not only minimizes the total unit cost but also achieves the most balanced cost structure, effectively controlling both electricity expenditure and membrane degradation. This demonstrates the value of integrating high-fidelity physical models for operational guidance, ensuring both economic efficiency and component longevity.
}
\end{enumerate}

\begin{figure}[!h]
 \centering
    \subfigure[Operating current density and temperature vary different strategies in the Houston electricity market. The white scatters indicate the low frequency operation point, and the dark scatters indicate the high frequency operation point. ]{
    \includegraphics[width=0.8\textwidth]{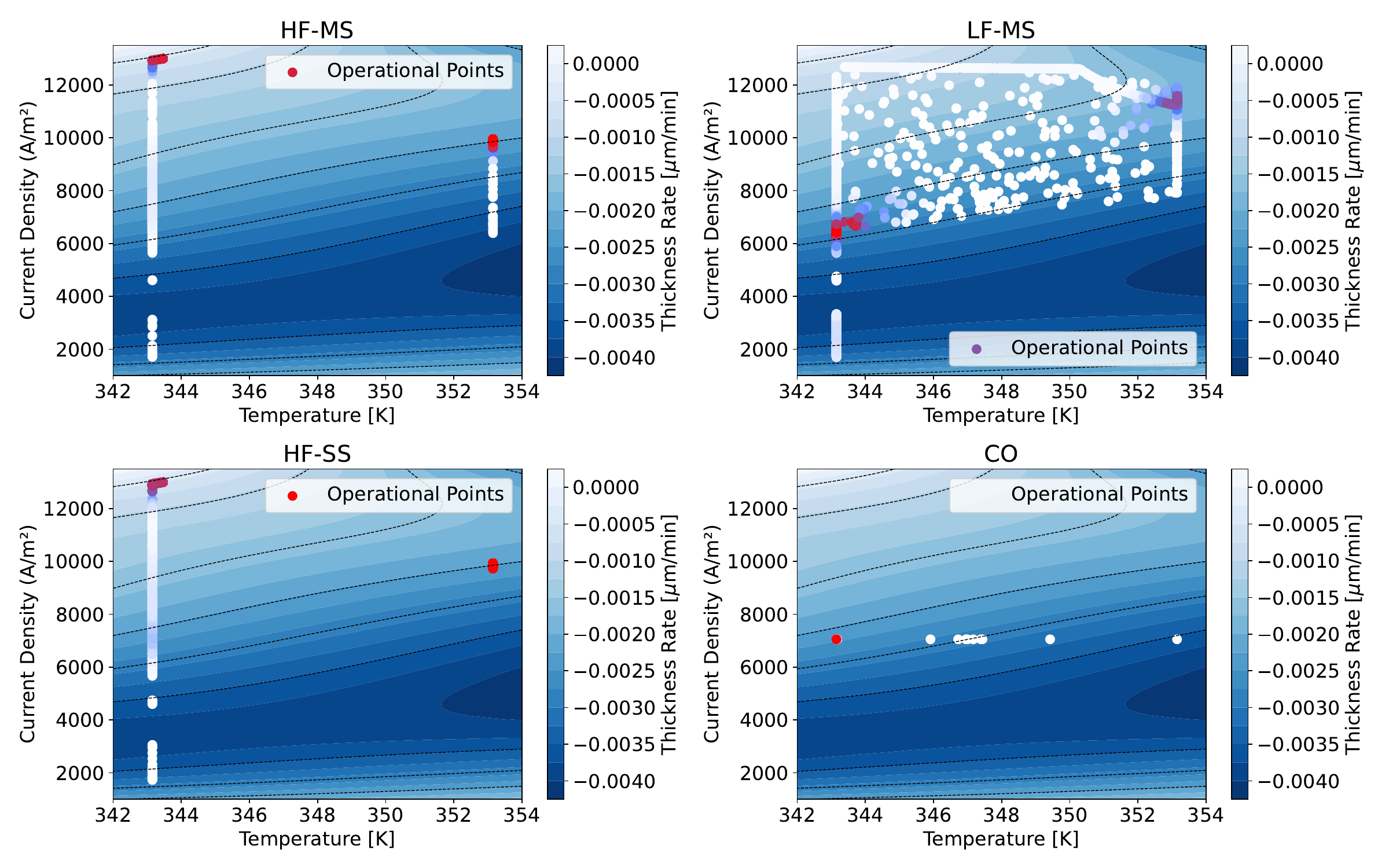} \label{fig:market_cost_breakdown1}
    } 
    \subfigure[The KDE of current density at 343.15 K, the black line is membrane degradation rate, and the red line is the KDE estimation of operating condition.]{
    \includegraphics[width=1.0\textwidth]{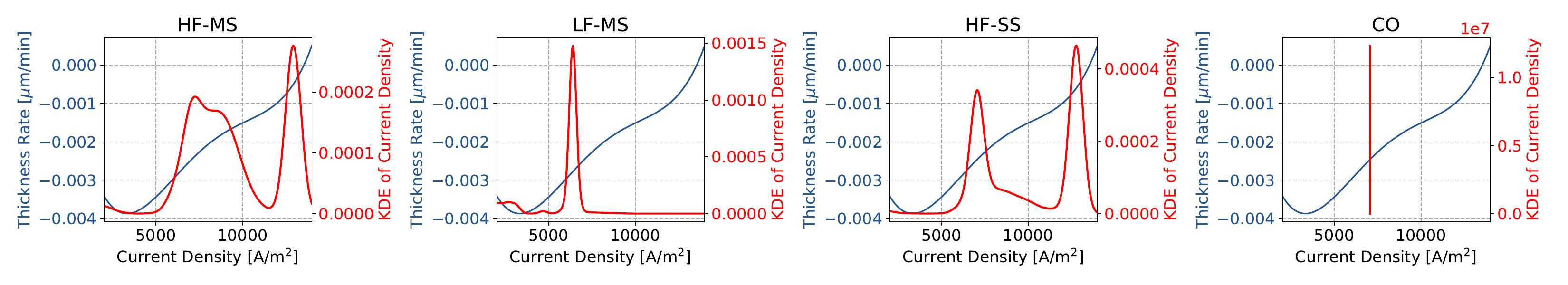} \label{fig:market_cost_breakdown2}
    } 
    \subfigure[The KDE of current density at 353.15 K, the black line is membrane degradation rate, and the red line is the KDE estimation of operating condition.]{
    \includegraphics[width=1.0\textwidth]{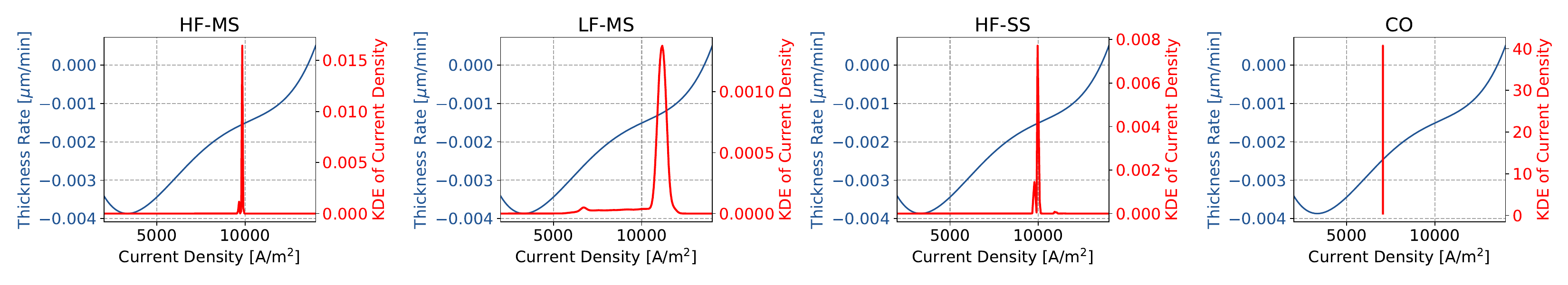} \label{fig:market_cost_breakdown3}
    } 
    \caption{Membrane degradation rate and operating temperature.}
    \label{fig:market_cost_breakdown}
\end{figure}
Finally, recognizing the significant contribution of membrane degradation to overall cumulative costs—particularly in the CO and HF-SS strategies (as shown previously)—we investigated the operational factors driving these cost disparities. Our goal was to elucidate how the HF-MS and HF-SS strategies effectively mitigate membrane degradation. To this end, we analyzed the operating conditions of each strategy in relation to the membrane degradation rate, as shown in~\Cref{fig:market_cost_breakdown1}. From~\Cref{fig:market_cost_breakdown1}, it is evident that both HF-MS and HF-SS consistently restrict operation to regions characterized by lower temperatures and moderate current densities, thereby avoiding regimes that accelerate membrane degradation. In contrast, LF-MS and CO display broader, more dispersed operational profiles, frequently operating in high-risk zones for membrane failure. 

Given that operation is concentrated at 343.15 K and 353.15 K, we further examined the distribution of current densities and corresponding membrane degradation rates at these temperatures in~\Cref{fig:market_cost_breakdown2,fig:market_cost_breakdown3}. Kernel density estimation (KDE) at 343.15 K (\Cref{fig:market_cost_breakdown2}) confirms that HF-MS and HF-SS strategies cluster around current densities associated with minimal membrane degradation, resulting in a clear advantage over LF-MS, with degradation rate differences on the order of 0.001 $\mu$m/min. These findings demonstrate that high-fidelity operation strategy—especially when combined with multi-scale market participation—strategically minimizes exposure to conditions that accelerate membrane failure, directly contributing to the observed reductions in cumulative membrane costs. At 353.15 K, even though, this advantage is less pronounced: LF-MS exhibits slightly lower degradation rates than the HF strategies, but the difference is marginal (less than 0.0005 $\mu$m/min). This suggests that while higher operating temperatures modestly increase degradation rates for HF strategies, disciplined current density management ensures that the overall impact remains small and comparable to LF strategies. Thus, the operational benefits of HF-MS and HF-SS strategies are largely preserved even at 353.15 K.
\section{Discussions}\label{sec:discussionsResults}
Recent advances in electricity market participation strategies and electrolyzer technologies have opened new opportunities for enhancing the competitiveness of green hydrogen. Our findings demonstrate that integrating high-fidelity electrolyzer modeling with dynamic, multi-scale market participation enables not only substantial reductions in electricity expenditure but also strategic minimization of membrane degradation—two key factors that have historically constrained the economic viability of PEM-based green hydrogen. By quantifying how optimized operational strategies (e.g., HF-MS) shift plant operation toward low-risk regimes—characterized by moderate current densities and lower temperatures, which reduce membrane degradation rates—we show that such strategies directly translate into a lower levelized cost of hydrogen compared to both traditional constant operation and lower-fidelity flexible operation approaches.

While our study highlights that disciplined operation—specifically, targeting moderate current densities and lower operating temperatures—can greatly extend membrane degradation cost without sacrificing market responsiveness, important challenges remain. Future work should focus on incorporating more accurate modeling of external wind and solar power availability to minimize reliance on electricity purchases from external markets~\cite{wang2018greening,ramakrishnan2024offshore,collins2025levelized}. Additionally, integrating electrolyzer operation with downstream chemical production processes~\cite{li2022bifunctional,ma2023exploiting,wei2024comparative} offers further opportunities for holistic system optimization and value creation.

The insights and modeling approaches developed here are broadly transferable to other regions and power market settings. By quantifying the tangible operational and economic benefits of coupling high-fidelity models with market-driven flexibility, our work provides a roadmap for hydrogen plant operators, technology developers, and market designers. Such strategies will be key to scaling green hydrogen and realizing its role as a cornerstone of decarbonized energy systems worldwide.

\section{Methods}

\subsection{Rolling Horizon Optimization Approach}
The operational strategy for the PEM electrolyzer plant in our work was determined using a rolling horizon optimization framework, which co-optimizes market participation and plant operation in response to dynamic conditions. The framework integrates a predictive controller with a dynamic plant simulator in a closed-loop configuration.

The predictive controller solves a finite-horizon optimization problem at each decision interval to determine the optimal control policy. The problem is formulated to minimize total operational costs, which include electricity procurement across both Day-Ahead and Real-Time Markets, as well as the projected cost of membrane degradation. The optimization is subject to a set of constraints that represent the physical and operational limits of the electrolyzer system. The general form of the optimization is expressed as:

\begin{equation}
\begin{aligned}
&\min_{x,u} ~ f(x, u) \\
&\text{subject to:} \\
&\qquad g(x, u) \leq 0 \hspace{2em} \text{(PEM Operational and Electricity Market Constraints)} \\
&\qquad h(x, u) = 0 \hspace{2.05em} \text{(Mass and Energy Balances of the PEM Plant)} \\
&\qquad \frac{dx}{dt} = q(x, u) \hspace{1.2em} \text{(System Dynamics)} \\
\end{aligned}
\end{equation}

Here, $x$ represents the set of state variables (e.g., membrane thickness, hydrogen storage level), and $u$ denotes the control actions (e.g., market bids, electricity strength). On this basis, the functions $g(\cdot)$, $h(\cdot)$, and $q(\cdot)$ define the system's operational envelope, conservation laws, and dynamic evolution, respectively. Notably, the $q(\cdot)$ term flexibly accommodates different model fidelities—including high-fidelity and low-fidelity models—enabling systematic comparison of operational and economic outcomes across multiple control and modeling paradigms. 

Based on this, the plant simulator module propagates the system's state forward in time based on the control policy determined by the controller. It uses a high-fidelity model that captures the plant's core physicochemical dynamics, including the electrochemical reactions, energy balances, and the empirical membrane degradation model. The system's differential equations are solved numerically using the forward Euler scheme~\cite{butcher2016numerical}.

This rolling horizon process is executed iteratively. At each step, the controller receives the current plant state from the simulator and the latest electricity prices from both the DAM and RTM. The optimization is performed over a horizon extending to 23:59 of the current day, generating an optimal sequence of control actions, including RTM bidding strategies. At 9:00 each day, the framework additionally determines the DAM bid quantities for the following day, in accordance with market requirements. At all other time steps, only RTM bidding and plant operation decisions are updated. Only the first action in this sequence is implemented. The simulator then calculates the resulting plant state for the next time step, which is fed back to the controller, and the cycle repeats. This closed-loop structure allows the strategy to adapt to evolving market prices and plant conditions. The optimization problem that determines the control policy and the simulator problem that determines the system response were solved using the IPOPT solver~\cite{wachter2006implementation} within the Pyomo optimization environment~\cite{hart2011pyomo,hart2011pyomo}.  All detailed model equations, parameters, and numerical solver settings are provided in the Appendix.

\subsection{Techno-Economic Data}
The techno-economic parameters in this study are primarily adapted from published literature and local market reports, with particular focus on the electricity markets of Houston and California. Historical DAM and RTM price data for 2022 were obtained from the Electric Reliability Council of Texas for Houston and the California Independent System Operator for California, capturing representative temporal price fluctuations for each region and ensuring robust analysis under typical market conditions.

PEM system modeling parameters, including electrolyzer stack performance, degradation rates, efficiency, and dynamic operational constraints, are based on reference~\cite{gorgun2006dynamic}. The modeled PEM plant is assumed to be utility-scale, with a nominal production capacity of 1 ton per hour. All capital, operation, maintenance, and electricity costs are annualized or accounted for, and levelized as the total hydrogen production cost per kilogram. 

To reflect actual market operation, our model incorporates both DAM and RTM electricity prices, with DAM cleared hourly and RTM cleared every 15 minutes. Each day, the PEM plant submits a 24-hour electricity consumption bid to the DAM for the next day at 9:00, and it is assumed that the scheduled electricity is fully delivered unless deviations occur in the RTM. Intra-hour adjustments are allowed based on 15-minute RTM prices, enabling operational flexibility while respecting DAM commitments and system constraints. All relevant grid transaction fees, market penalties, and regulatory charges in Houston and California are included in the cost model.

All techno-economic assumptions, input parameters, and Python modeling scripts necessary for reproducing the results are provided in the Appendix. Detailed PEM system specifications and a breakdown of cost assumptions are described in Table A.1 and Section A of the Appendix.

\bibliographystyle{unsrtnat}
\nobibliography*
\bibliography{ref_adjoint}
\newpage
\appendix
\counterwithin{equation}{section}
\setcounter{equation}{0}
\counterwithin{table}{section}
\setcounter{table}{0}
\counterwithin{figure}{section}
\setcounter{figure}{0}
\begin{appendices}
    \let\oldsection\section
    \renewcommand{\section}[1]{
        \oldsection{\textcolor{black}{#1}}
        \addcontentsline{apc}{section}{\protect\numberline{\thesection}\textcolor{blue}{#1}}
    }

    \let\oldsubsection\subsection
    \renewcommand{\subsection}[1]{
        \oldsubsection{\textcolor{black}{#1}}
        \addcontentsline{apc}{subsection}{\protect\numberline{\thesubsection}\textcolor{blue}{#1}}
    }



\section{Additional Controller and Simulator Information}
 Detailed information on the plant and controller models for the PEM-based electrolyzer is elaborated in this section. 

\subsection{Plant Model}\label{sec:plantModel}
\subsubsection{Mass Conservation of PEM Electrolyzer}\label{subsec:HydrogenProductionConservation}
In this section, we present the mass conservation equations for the PEM electrolyzer. Specifically, the generation rates of hydrogen at the cathode and oxygen at the anode are given by:
\begin{subequations}
\begin{align}\label{eq:cathodeEquation}\text{cathode:~~}  \dot{\mathrm{N}}_{\text{H}_2,\text{gen}}= \dfrac{nI\eta}{2\mathcal{F}},\end{align}
\begin{align}\label{eq:anodeEquation}\text{anode:~~}  \dot{\mathrm{N}}_{\text{O}_2,\text{gen}} = \dfrac{nI\eta}{4\mathcal{F}},\end{align}
\end{subequations}
where the subscript `gen' is the abbreviation of generation, $n$ is the stacked number of the PEM electrolyzer, $I$ is the operation current, $\eta$ is the efficiency of the electrolyzer, and $\mathcal{F}$ is the celebrated Fraday constant. Based on this, the net product flow rates at the cathode and anode are then expressed as:
\begin{subequations}
\begin{align}\label{eq:cathodeMassEquation}\text{cathode:~~}  \dot{\mathrm{N}}_{\text{H}_2,\text{net}}=    \dot{\mathrm{N}}_{\text{H}_2,\text{gen}} -   \dot{\mathrm{N}}_{\text{H}_2,\text{out}},\end{align}
\begin{align}\label{eq:anodeMassEquation}\text{anode:~~}   \dot{\mathrm{N}}_{\text{O}_2,\text{net}}=    \dot{\mathrm{N}}_{\text{O}_2,\text{gen}} -   \dot{\mathrm{N}}_{\text{O}_2,\text{out}}.\end{align}
\end{subequations}
Notably, the produced hydrogen can be distributed to two destinations: the hydrogen storage tank and the downstream industrial plant. This allocation is represented by:
\begin{equation}\label{eq:chamberOutEquation}
\dot{\mathrm{N}}_{\text{H}_2,\text{out}} = \dot{\mathrm{N}}_{\text{H}_2,\text{stor,in}} +  \dot{\mathrm{N}}_{\text{H}_2,\text{el,plant}}.
\end{equation}
The net flow rate into the hydrogen storage tank is given by:
\begin{equation}\label{eq:sotrageEquation}
 \dot{\mathrm{N}}_{\text{H}_2,\text{stor,in}} -  \dot{\mathrm{N}}_{\text{H}_2,\text{stor,out}} =  \dot{\mathrm{N}}_{\text{H}_2,\text{stor,net}}
\end{equation}
Accordingly, the total hydrogen supplied to the downstream plant is:
\begin{equation}\label{eq:plantReceiveEquation}
    \dot{\mathrm{N}}_{\text{H}_2,\text{plant}} = \dot{\mathrm{N}}_{\text{H}_2,\text{el,plant}} + \dot{\mathrm{N}}_{\text{H}_2,\text{stor,out}}.
\end{equation}

To ensure a stable supply of hydrogen to the downstream process, the following equality constraint is imposed in the model:
\begin{equation}\label{eq:setpointEquation}
     \dot{\mathrm{N}}_{\text{H}_2,\text{plant}} =  \dot{\mathrm{N}}_{\text{H}_2,\text{setpoint}}
\end{equation}
Similarly, the hydrogen storage tank is subject to volume constraints, which are represented by the following inequality:
\begin{equation}\label{ineq:1}
  {\mathrm{N}}_{\text{H}_2,\text{stor,lb}} \le  {\mathrm{N}}_{\text{H}_2,\text{stor,net} }
\le  {\mathrm{N}}_{\text{H}_2,\text{stor,ub}},
\end{equation}
where $ {\mathrm{N}}_{\text{H}_2,\text{stor,lb}}$ and $ {\mathrm{N}}_{\text{H}_2,\text{stor,ub}}$ denote the lower and upper bounds of the hydrogen storage tank, respectively. In addition, to ensure the safe operation of the PEM electrolyzer, the hydrogen production rate is constrained within its minimum and maximum limits:
\begin{equation}\label{ineq:2}
     \dot{\mathrm{N}}_{\text{H}_2,\text{gen}}^{\text{min}}\le      \dot{\mathrm{N}}_{\text{H}_2,\text{gen}}   \le       \dot{\mathrm{N}}_{\text{H}_2,\text{gen}}^{\text{max}}.
\end{equation}
With these equations, mass conservation constraints for the PEM plant are now established.

\subsubsection{Energy Conservation of the Whole Plant}
While \Cref{subsec:HydrogenProductionConservation} effectively models the mass conservation of the PEM plant, it does not explicitly account for energy conservation. To address this gap, the energy conservation is further elaborated in this part. First, we present the overall voltage for a single PEM stack, as outlined in reference \cite{XIONG2024251}:
\begin{equation}\label{eq:totalVoltage}
    V_{\text{total}} = V_{\text{act} } +  V_{\text{oc} } +  V_{\text{ohm} },
\end{equation}
where $V$ represents voltage, $V_{\text{total}}$ is the total voltage, $V_{\text{act}}$ is the activation voltage, $V_{\text{oc}}$ is the open-circuit voltage, and $V_{\text{ohm}}$ is the Ohmic voltage. Each term in \Cref{eq:totalVoltage} can be defined as follows:
\begin{itemize}
    \item{The expression for $V_{\text{act} }$ is given as follows:
\begin{equation}\label{eq:actVoltage}
    V_{\text{act} } = \frac{\mathcal{R}\mathcal{T}}{2\mathcal{F}\mathcal{C}_{\text{charge}}}\ln{\frac{I}{\rho_I\mathcal{A}}},
\end{equation}
where $\mathcal{R}$ is the celebrated Avogadro constant, $\mathcal{T}$ is the operating temperature of the electrolyzer, $\mathcal{C}_{\text{charge}}$ is the charge coefficient, $\mathcal{A}$ is the membrane area, $\rho_I$ is the exchange current density. }
\item{The expression for $V_{\text{oc} }$ is given as follows: 
\begin{equation}\label{eq:ocVoltage}
       V_{\text{oc} } =  \dfrac{R\mathcal{T}}{2\mathcal{F}}\ln{p_{\text{H}_2} \sqrt{p_{\text{O}_2}} }+V_{\text{rev} } , 
\end{equation}
where $p_{\text{H}_2}$ and $p_{\text{O}_2}$ are the pressure of hydrogen and oxygen, which are computed by the following equations with the help of electrolyzer volume $\mathcal{V}$: 
\begin{subequations}
\begin{align}\label{eq:h2Pressure}
    p_{\text{H}_2} = \dfrac{\mathrm{N}_{\text{H}_2,\text{net}}\mathcal{R}\mathcal{T}}{\mathcal{V}} ,
\end{align}
\begin{align}\label{eq:o2Pressure}
     p_{\text{O}_2} = \dfrac{\mathrm{N}_{\text{O}_2,\text{net}}\mathcal{R}\mathcal{T}}{\mathcal{V}},
\end{align}
\end{subequations}
and $V_{\text{rev} }$ is the reversible potential, which is delineated by the following equation:
\begin{equation}\label{eq:revVoltage}
  V_{\text{rev} } =  1.299 - 0.9 \times 10^{-3} \times  (\mathcal{T} - 298) .
\end{equation}
}
\item{The expression for $V_{\text{ohm} }$ is given as follows: 
\begin{equation}\label{eq:ohmVoltage}
       V_{\text{ohm} } = \frac{I\times\varepsilon}{\mathcal{A}\times\beta}, 
\end{equation}
where $\varepsilon$ is the membrane thickness, and $\beta$ is the membrane conductivity, which satisfies the following equation:
\begin{equation}\label{eq:membraneConductivity}
    \beta = (0.00514\lambda_{\text{mem}}-0.00326)\exp\left[1268\left(\dfrac{1}{303}-\dfrac{1}{\mathcal{T}}\right)\right],
\end{equation}
and $\lambda_\text{mem}$ is the water content of PEM.

}
\end{itemize}
In practical applications, it is crucial for an electrolyzer to operate within specific voltage limits to ensure safe and efficient performance, for the sake of helping prevent damage to the system and enhancing the longevity of the equipment. Therefore, the following inequality is introduced to restrict the operation condition of $V_{\text{total}}$:
\begin{equation}\label{eq:totalVoltageBound}
 V_{\text{total}}^{\text{min}} \le   V_{\text{total}} \le  V_{\text{total}}^{\text{max}}.
\end{equation}
where $V_{\text{total}}^{\text{min}}$ represents the minimum operational voltage, below which the electrolyzer may not function effectively, leading to insufficient hydrogen production and potential system instability. Conversely, $V_{\text{total}}^{\text{max}}$ denotes the maximum operational voltage, which, if exceeded, could result in excessive heating, increased wear on the components, or even catastrophic failure.

Building on this foundation, the total power consumption of the PEM plant can be expressed by the following equation:
\begin{equation}\label{eq:plantPower}
    P_{\text{plant}} =   P_{\text{PEM}} +  P_{\text{extra}}.
\end{equation}
In this equation, $P_{\text{PEM}}$ represents the energy consumption of the PEM electrolyzer, which is calculated using the following formula:
\begin{equation}\label{eq:PEMtotalVoltageCompute}
    P_{\text{PEM}} = V_{\text{total}} \times I\times n,
\end{equation}
where it is assumed that the PEM electrolyzers are connected in series. In addition, the term $P_{\text{extra}}$ accounts for the additional energy required to maintain the operational conditions necessary for hydrogen production. This includes the power consumed by auxiliary devices such as heat exchangers, compressors, and cooling systems, etc. In this manuscript, $P_{\text{extra}}$ is estimated as follows:
\begin{equation}\label{eq:extraPlantConsumption}
    P_{\text{extra}} = \mathcal{K}_\text{extra} \dfrac{\mathrm{d} \mathrm{N}_{\text{H}_2,\text{gen}}}{\mathrm{d}t},
\end{equation}
where ${P}_\text{extra}$ is a proportionality constant that indicates the energy required per unit of hydrogen produced.

\subsubsection{Membrane Degradation Model}
In practice, the performance of PEM electrolyzers inevitably deteriorates over time due to various degradation mechanisms. Among these, membrane thinning has been identified in the literature as one of the most prominent manifestations of degradation\footnote{As noted in reference~\cite{chandesris2015membrane}, membrane degradation does not affect the membrane area; thus, the influence of membrane area is explicitly omitted in this manuscript.}. Membrane thinning poses significant safety and economic challenges, including an increased risk of hazardous gas crossover (potentially leading to explosions) and higher costs associated with membrane replacement. To systematically address these concerns, this subsection focuses on modeling the evolution of membrane thickness during plant operation.

To the best of our knowledge, a rigorous mathematical model describing membrane thickness evolution is currently lacking in the literature. Therefore, we propose to develop an empirical model based on available references. Specifically, following the approach in reference~\cite{chandesris2015membrane}, we model the membrane thinning rate, denoted as $\dot{\varepsilon}$, as a function of operating temperature $\mathcal{T}$ and current density $\frac{I}{\mathcal{A}}$. According to the experimental data reported in~\cite{chandesris2015membrane}, the membrane thinning rate $\dot{\varepsilon}$ can be expressed as follows:
\begin{equation}\label{eq:membraneDegradationEquationResult}
\begin{aligned}
        \dot{\varepsilon} =& (-0.008255 \times \mathcal{T} + 2.906615) \times (\dfrac{I}{\mathcal{A}} )^4
            + (0.021855 \times \mathcal{T} - 7.740815) \times (\dfrac{I}{\mathcal{A}} ) ^3\\
            & + (-0.01798 \times \mathcal{T} + 6.44534) \times (\dfrac{I}{\mathcal{A}} ) ^2
            + (0.00415 \times \mathcal{T} - 1.53825) \times (\dfrac{I}{\mathcal{A}} ) ^1\\
    &        + (-0.00005 \times \mathcal{T} + 0.01715).
\end{aligned}
\end{equation}


\subsection{Controller Model}\label{sec:controllerModelResult}
In addition to the PEM model presented in~\Cref{sec:plantModel}, the controller model involved in this paper also includes the constraints on the electricity market participant and the related cost. 
\subsubsection{Electricity Market Bidding}
to fully utilize the potential of the both electricity market and realize the goal of money saving, we first establish the energy conservation from these two energy markets as follows:
\begin{equation}\label{eq:DAMRTMConservation}
    P_{\text{DAM},t} + P_{\text{RTM},t} = P_{\text{plant},t}\times9\times10^{-4},
\end{equation}
where $9\times10^{-4}$ serves as the energy conversion coefficient, translating the power unit (W) into energy consumption (MW$\cdot$hr).

In this procedure, to prevent the downtime of the PEM plant, the following inequality constraint is introduced to maintain the power consumption of the PEM at a level exceeding 10\% of its maximum capacity while not exceeding the maximum power limit:
\begin{equation}\label{eq;plantBound}
   10\%\times(0.25\times P_{\text{plant}}^{\text{max}})\le P_{\text{plant},t} \le 0.25\times P_{\text{plant}}^{\text{max}},
\end{equation}
where $0.25$ acts as the energy conversion coefficient, converting hourly energy consumption into quarterly energy consumption in units of MW$\cdot$hr, $P_{\text{plant}}^{\text{max}}$ is the maximum capacity of the PEM plant, which can be estimated by the following equation: 
\begin{equation}\label{eq:maximumConsumption}
      P_{\text{PEM},\max} = \mathcal{K}_{\text{H}_2}\times  \dot{ \mathrm{N}}_{\text{H}_2,\text{gen}}^\text{max},
\end{equation}
and $\mathcal{K}_{\text{H}_2}$ represents the energy consumption per mole of hydrogen generated. On this basis, the following inequality constraint is introduced to prohibit the sale of electricity to the DAM and to prevent excessive purchases from the DAM:
\begin{equation}\label{eq:DAMBound}
  0\le  P_{\text{DAM},t} \le P_{\text{plant}}^{\text{max}}.
\end{equation}
Similarly, the following inequality constraint is established to prevent the overselling of electricity to the RTM and to avoid excessive purchases from the RTM:
\begin{equation}\label{eq:RTMBound}
  -90\%\times P_{\text{plant}}^{\text{max}}\le  P_{\text{DAM},t} \le P_{\text{plant}}^{\text{max}}.
\end{equation}

\subsubsection{
Cost Computation}
The cost for the whole plant consists of the electricity cost and the PEM component cost. Specifically, the electricity cost of the PEM plant along the whole control horizon $\mathrm{T}$ can be given as follows:
\begin{equation}\label{eq:totalElecCost}
    C_{\text{elec}} =  
    \sum_{t=0}^{\mathrm{T}-1}{\left[C_{\text{RTM},t} (P_{\text{RTM},t}) + C_{\text{DAM},t} (P_{\text{DAM},t})\right]},
\end{equation}
where $ C_{\text{DAM},t}$ and $C_{\text{RTM},t}$ are the electricity price of DAM and RTM at moment $t$, respectively. Meanwhile, the membrane degradation cost can be given as follows:
\begin{equation}\label{eq:totalMembraneCost}
    C_{\text{mem,HF}} = nP_{\text{mem}} (\varepsilon_{\mathrm{T}} - \varepsilon_{\mathrm{0}} ),
\end{equation}
where $ P_{\text{mem}}$ is the membrane cost coefficient. 

Based on the preceding contents, the optimal control problem is formulated for the control policy determination. Specifically, the objective is to identify the optimal control policy $u_t$ over the planning horizon so as to minimize the total operational cost, which includes membrane replacement costs ($C_{\text{mem}}$) and electricity procurement costs ($C_{\text{elec}}$). In the context of this study, the control policy $u_t$ corresponds to the operational condition of the PEM electrolyzer at each time step. Based on the abovementioned contents, we therefore formulate the programming problem for the HF-MS operation strategy can be given as follows:
\begin{equation}\label{eq:OptimalControlProblem2}
\begin{footnotesize}
\text{HF-MS:~~}
\begin{aligned}
  &\mathop{\arg\min}_{{u}_t}\quad  
 \text{Obj}= C_{\text{mem}} + C_{\text{elec}} \\
  & \text{s.t.}  \begin{cases}
  & {u}_t =\left\{  P_{\text{DAM},t},P_{\text{RTM},t},\mathcal{T}_t,\mathcal{I}_t,\dot{\mathrm{N}}_{\text{H}_2,\text{el,plant}},\dot{\mathrm{N}}_{\text{H}_2,\text{stor,in}},\dot{\mathrm{N}}_{\text{H}_2,\text{stor,out}}\right\},\\
  &\text{Mass Conservation:}~\text{\Cref{eq:cathodeEquation,eq:anodeEquation,eq:cathodeMassEquation,eq:anodeMassEquation,eq:chamberOutEquation,eq:sotrageEquation,eq:plantReceiveEquation,eq:setpointEquation}} \\
  &\text{Energy Conservation:}~\text{\Cref{eq:totalVoltage,eq:actVoltage,eq:ocVoltage,eq:h2Pressure,eq:o2Pressure,eq:revVoltage,eq:ohmVoltage,eq:membraneConductivity,eq:plantPower,eq:PEMtotalVoltageCompute,eq:extraPlantConsumption}} \\
  &\text{Membrane Degradation:}~\text{\Cref{eq:membraneDegradationEquationResult}}\\
  &\text{Electricity Market Trading:}~\text{\Cref{eq:DAMRTMConservation}}\\
&\text{Cost Computation:}~\text{\Cref{eq:totalElecCost,eq:totalMembraneCost}}\\
   &\text{Inequality Constraints:}~\text{\Cref{ineq:1,ineq:2,eq:totalVoltageBound,eq;plantBound,eq:DAMBound,eq:RTMBound}}
    \end{cases}
\end{aligned}
\end{footnotesize}
\end{equation}
Similarly, compared with the HF-MS, the HF-SS merely requires the following constraint:
\begin{equation}\label{eq:freeRTM}
    P_{\text{RTM},t} = 0,
\end{equation}
and the corresponding programming problem can be given as follows:
\begin{equation}\label{eq:OptimalControlProblem}
\begin{footnotesize}
\text{HF-SS:~~}
\begin{aligned}
  &\mathop{\arg\min}_{{u}_t}\quad  
  \text{Obj}=C_{\text{mem}} + C_{\text{elec}} \\
  & \text{s.t.}  \begin{cases}
  & {u}_t =\left\{  P_{\text{DAM},t},P_{\text{RTM},t},\mathcal{T}_t,\mathcal{I}_t,\dot{\mathrm{N}}_{\text{H}_2,\text{el,plant}},\dot{\mathrm{N}}_{\text{H}_2,\text{stor,in}},\dot{\mathrm{N}}_{\text{H}_2,\text{stor,out}}\right\},\\
  &\text{Mass Conservation:}~\text{\Cref{eq:cathodeEquation,eq:anodeEquation,eq:cathodeMassEquation,eq:anodeMassEquation,eq:chamberOutEquation,eq:sotrageEquation,eq:plantReceiveEquation,eq:setpointEquation}} \\
  &\text{Energy Conservation:}~\text{\Cref{eq:totalVoltage,eq:actVoltage,eq:ocVoltage,eq:h2Pressure,eq:o2Pressure,eq:revVoltage,eq:ohmVoltage,eq:membraneConductivity,eq:plantPower,eq:PEMtotalVoltageCompute,eq:extraPlantConsumption,eq:freeRTM}} \\
  &\text{Membrane Degradation:}~\text{\Cref{eq:membraneDegradationEquationResult}}\\
  &\text{Electricity Market Trading:}~\text{\Cref{eq:DAMRTMConservation}}\\
&\text{Cost Computation:}~\text{\Cref{eq:totalElecCost,eq:totalMembraneCost}}\\
   &\text{Inequality Constraints:}~\text{\Cref{ineq:1,ineq:2,eq:totalVoltageBound,eq;plantBound,eq:DAMBound,eq:RTMBound}}
    \end{cases}
\end{aligned}
\end{footnotesize}
\end{equation}

Notably, for LF-MS operation strategy, we regardless the membrane degradation equation and treat the membrane cost as the coefficient of the hydrogen production flowrate as follows:
\begin{equation}\label{eq:lfMembraneCost}
     C_{\text{mem,LF}} =\mathcal{K}_{\text{mem,H}_2}\sum_{t=0}^{\mathrm{T}-1}{\dot{\mathrm{N}}_{\text{H}_2,\text{gen}}},
\end{equation}
where $\mathcal{K}_{\text{mem,H}_2}$ is the hydrogen flowrate coefficient. Consequently, the programming problem for the LF-MS can be given as follows:
\begin{equation}\label{eq:OptimalControlProblem4}
\begin{footnotesize}
\text{LF-MS:~~}
\begin{aligned}
  &\mathop{\arg\min}_{{u}_t}\quad  \text{Obj}=  C_{\text{mem,LF}}+ C_{\text{elec}} \\
   & \text{s.t.}  \begin{cases}
  & {u}_t =\left\{  P_{\text{DAM},t},P_{\text{RTM},t},\mathcal{T}_t,\mathcal{I}_t,\dot{\mathrm{N}}_{\text{H}_2,\text{el,plant}},\dot{\mathrm{N}}_{\text{H}_2,\text{stor,in}},\dot{\mathrm{N}}_{\text{H}_2,\text{stor,out}}\right\},\\
  &\text{Mass Conservation:}~\text{\Cref{eq:cathodeEquation,eq:anodeEquation,eq:cathodeMassEquation,eq:anodeMassEquation,eq:chamberOutEquation,eq:sotrageEquation,eq:plantReceiveEquation,eq:setpointEquation}} \\
  &\text{Energy Conservation:}~\text{\Cref{eq:totalVoltage,eq:actVoltage,eq:ocVoltage,eq:h2Pressure,eq:o2Pressure,eq:revVoltage,eq:ohmVoltage,eq:membraneConductivity,eq:plantPower,eq:PEMtotalVoltageCompute,eq:extraPlantConsumption}} \\
  &\text{Electricity Market Trading:}~\text{\Cref{eq:DAMRTMConservation}}\\
&\text{Cost Computation:}~\text{\Cref{eq:totalElecCost,eq:lfMembraneCost}}\\
   &\text{Inequality Constraints:}~\text{\Cref{ineq:1,ineq:2,eq:totalVoltageBound,eq;plantBound,eq:DAMBound,eq:RTMBound}}
    \end{cases}
\end{aligned}
\end{footnotesize}
\end{equation}
For the constant open problem, we fix part of the control policy as follows:
\begin{equation}\label{eq:constOpen}
\begin{cases}
    \mathcal{I}_t &=  3.526\times 10^4\\
    \dot{\mathrm{N}}_{\text{H}_2,\text{el,plant}} &=  0\\
 \dot{\mathrm{N}}_{\text{H}_2,\text{stor,out}}& =  0\\
\end{cases},
\end{equation}
and formulate the following optimization as follows:
\begin{equation}\label{eq:OptimalControlProblem3}
\begin{footnotesize}
\text{CO:~~}
\begin{aligned}
  &\mathop{\arg\min}_{{u}_t}\quad   \text{Obj}=C_{\text{mem,LF}}+ C_{\text{elec}} \\
   & \text{s.t.}  \begin{cases}
  & {u}_t =\left\{  P_{\text{DAM},t},P_{\text{RTM},t},\mathcal{T}_t,\mathcal{I}_t,\dot{\mathrm{N}}_{\text{H}_2,\text{el,plant}},\dot{\mathrm{N}}_{\text{H}_2,\text{stor,in}},\dot{\mathrm{N}}_{\text{H}_2,\text{stor,out}}\right\},\\
  &\text{Mass Conservation:}~\text{\Cref{eq:cathodeEquation,eq:anodeEquation,eq:cathodeMassEquation,eq:anodeMassEquation,eq:chamberOutEquation,eq:sotrageEquation,eq:plantReceiveEquation,eq:setpointEquation}} \\
  &\text{Energy Conservation:}~\text{\Cref{eq:totalVoltage,eq:actVoltage,eq:ocVoltage,eq:h2Pressure,eq:o2Pressure,eq:revVoltage,eq:ohmVoltage,eq:membraneConductivity,eq:plantPower,eq:PEMtotalVoltageCompute,eq:extraPlantConsumption}} \\
    &\text{Constant Open Constraint:}~\text{\Cref{eq:constOpen}} \\
  &\text{Electricity Market Trading:}~\text{\Cref{eq:DAMRTMConservation}}\\
&\text{Cost Computation:}~\text{\Cref{eq:totalElecCost,eq:lfMembraneCost}}\\
   &\text{Inequality Constraints:}~\text{\Cref{ineq:1,ineq:2,eq:totalVoltageBound,eq;plantBound,eq:DAMBound,eq:RTMBound}}
    \end{cases}
\end{aligned}
\end{footnotesize}
\end{equation}
Based on the abovementioned contents, the parameters for optimization are listed in~\Cref{tab:baselineComparison}. Based on this, we employ the IpOpt solver~\cite{wachter2006implementation} to address the programming, where the system dynamics are discretized using the forward Euler's method~\cite{butcher2016numerical}.

\begin{table}[htbp]
\centering
\caption{Detailed Parameters of the PEM Plant Modeling}\label{tab:baselineComparison}
\begin{threeparttable}
\begin{tabular}{l|c|l|l}
\toprule
Parameter & Symbol& Value & Unit\\ \midrule
PEM Stack Number &$n$ & 800 & Number \\
Electrolyzer Efficiency &$\eta$& 95 & \% \\
Membrane Area &$\mathcal{A}$& 5 & m$^\text{2}$ \\
Membrane Initial Thickness & $\varepsilon_0$ & 178 & $\mu$m \\
Exchange Current Density & $\rho_I$ & 0.00001 &  A$\cdot \text{cm}^{-1}$ \\
Water Content of PEM & $\lambda_{\text{mem}}$ & 14 &  - \\
Minimum Voltage & $V_{\text{total}}^{\text{min}} $ & 1.4 & V \\
Maximum Voltage & $V_{\text{total}}^{\text{max}} $ & 2.8 & V \\
Maximum Storage Volume & $\mathrm{N}_{\text{H}_2,\text{ms}}$& 7000 & $\text{kmol}$ \\
Minimum Storage Capacity Percentage & $C_{\text{min}}$& 21\% & \% \\
Maximum Storage Capacity Percentage & $C_{\text{max}}$& 100\% & \% \\
Hydrogen Provide Setpoint & $\dot{\mathrm{N}}_{\text{H}_2,\text{setpoint}}$& 500 & $\text{kmol}\cdot\text{hr}^{\text{-1}}$ \\
Hydrogen Generation Lower Bound & $\dot{\mathrm{N}}_{\text{H}_2,\text{gen}}^\text{min}$ & 100 & $\text{kmol}\cdot\text{hr}^{\text{-1}}$\\
Hydrogen Generation Upper Bound & $\dot{\mathrm{N}}_{\text{H}_2,\text{gen}}^{\text{max}}$ & 1000 & $\text{kmol}\cdot\text{hr}^{\text{-1}}$ \\
Current Density Lower Bound &$\left( \dfrac{I}{\mathcal{A}}\right)^\text{min}$ & 1000  & $\text{A}\cdot \text{m}^{-2}$ \\
Current Density Upper Bound &$\left( \dfrac{I}{\mathcal{A}}\right)^\text{max}$ & 13000 &$\text{A}\cdot \text{m}^{-2}$  \\
Hydrogen Flowrate Coefficient &$\mathcal{K}_{\text{mem,H}_2}$ & 1.388 &$\$\cdot \text{kmol}\cdot \text{hr}^{-1}$  \\
Temperature Lower Bound &$\mathcal{T}^{\text{min}}$ & 343  & $\text{K}$ \\
Temperature Upper Bound &$\mathcal{T}^{\text{max}}$ & 353 & $\text{K}$ \\
Membrane Cost Coefficient & $P_{\text{mem}}$& 203142 & $\$\cdot\mu\text{m}^{-1}$  \\
Extra Energy Consumption Coefficient &$\mathcal{P}_\text{extra}$ & 10 & $\text{kW}\cdot\text{hr}\cdot\text{kg}^{-1}$ \\
Energy Consumption Coefficient &$\mathcal{K}_{\text{H}_2}$ & 10 & $\text{kW}\cdot\text{hr}\cdot\text{kmol}^{-1}$ \\
Maximum Power of Plant &$ P_{\text{plant}}^{\text{max}}$ & 110000 &$\text{kW}\cdot\text{hr}$  \\
DAM Price at moment $t$ & $P_{\text{DAM},t}$&- & $\$\cdot \left(\text{MW}\cdot \text{hr}\right)^{-1}$ \\
RTM Price at moment $t$ & $P_{\text{RTM},t}$&- & $\$\cdot \left(\text{MW}\cdot \text{hr}\right)^{-1}$ \\
\bottomrule
\end{tabular}

  \end{threeparttable}
\end{table}

\subsection{Simulator Model}
In this context, the simulation problem can be represented as a special case of an optimization problem with zero degrees of freedom, meaning that all variables are fully determined by the constraints. In such cases, a dummy objective function $\text{Obj}$ is introduced—typically set to zero—solely to conform with standard optimization solver requirements, rather than to reflect any meaningful, real-world criterion. Accordingly, the simulator model can be formulated as follows:
\begin{equation}\label{eq:OptimalControlProblem1}
\begin{footnotesize}
\text{Simulator:~~}
\begin{aligned}
  &\mathop{\arg\min}_{{u}_t}\quad  
 \text{Obj} \equiv 0\\
  & \text{s.t.}  \begin{cases}
  & {u}_t =u_{t,\text{controller}},\\
  &\text{Mass Conservation:}~\text{\Cref{eq:cathodeEquation,eq:anodeEquation,eq:cathodeMassEquation,eq:anodeMassEquation,eq:chamberOutEquation,eq:sotrageEquation,eq:plantReceiveEquation,eq:setpointEquation}} \\
  &\text{Energy Conservation:}~\text{\Cref{eq:totalVoltage,eq:actVoltage,eq:ocVoltage,eq:h2Pressure,eq:o2Pressure,eq:revVoltage,eq:ohmVoltage,eq:membraneConductivity,eq:plantPower,eq:PEMtotalVoltageCompute,eq:extraPlantConsumption,eq:freeRTM}} \\
  &\text{Membrane Degradation:}~\text{\Cref{eq:membraneDegradationEquationResult}}\\
  &\text{Electricity Market Trading:}~\text{\Cref{eq:DAMRTMConservation}}\\
&\text{Cost Computation:}~\text{\Cref{eq:totalElecCost,eq:totalMembraneCost}}\\
    \end{cases}
\end{aligned}
\end{footnotesize}
\end{equation}
where $u_{t,\text{controller}}$ is the control policy obtained by the controller model. Based on this, we design the IpOpt solver~\cite{wachter2006implementation} to address the programming, where the system dynamics are discretized using the forward Euler's method~\cite{butcher2016numerical}. 

\subsection{Overall Workflow for Rolling Horizon Optimization}
\begin{figure}[!h]
 \centering
    \includegraphics[width=0.50\textwidth]{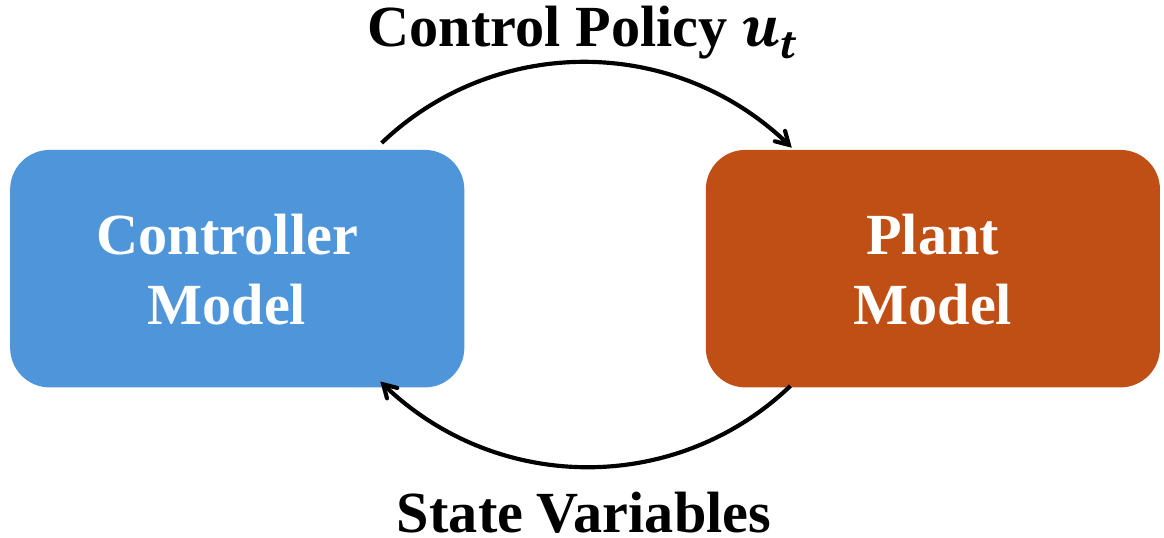}
        \caption{The illustration of the controller and plant.}
    \label{fig:overall_workflow_result}
\end{figure}

Based on the plant and controller models given in~\Cref{sec:plantModel,sec:controllerModelResult}, we further elaborate the workflow for the rolling horizon optimization. Specifically, the rolling horizon optimization is an interactive process between the controller and the model. The overall workflow of our methodology is illustrated in~\Cref{fig:overall_workflow_result}. This process is centered on an iterative interaction between a controller model and a plant model~\cite{rawlings2000tutorial}. Initially, the controller model, highlighted in the blue zone of the diagram, is solved to determine the optimal control policy, $u_t$. Subsequently, this determined control policy $u_t$ is applied to the plant model, which simulates the dynamic response of the PEM system highlighted in the oragne zone.The simulation captures the dynamic evolution of critical state variables, such as membrane thickness and hydrogen storage level, along with other essential parameters that delineate system performance. A key feature of this workflow is the feedback mechanism, where these updated state variables from the plant model are relayed back as inputs to the controller model for subsequent optimization cycles. Through the recursive execution of these steps—optimizing the control policy based on current states and then simulating the plant's response—the system progressively converges towards the optimal operational strategy for the entire plant.

\section{Additional Empirical
Experimental Results}
\begin{figure}[!h]
 \centering
    \includegraphics[width=0.95\textwidth]{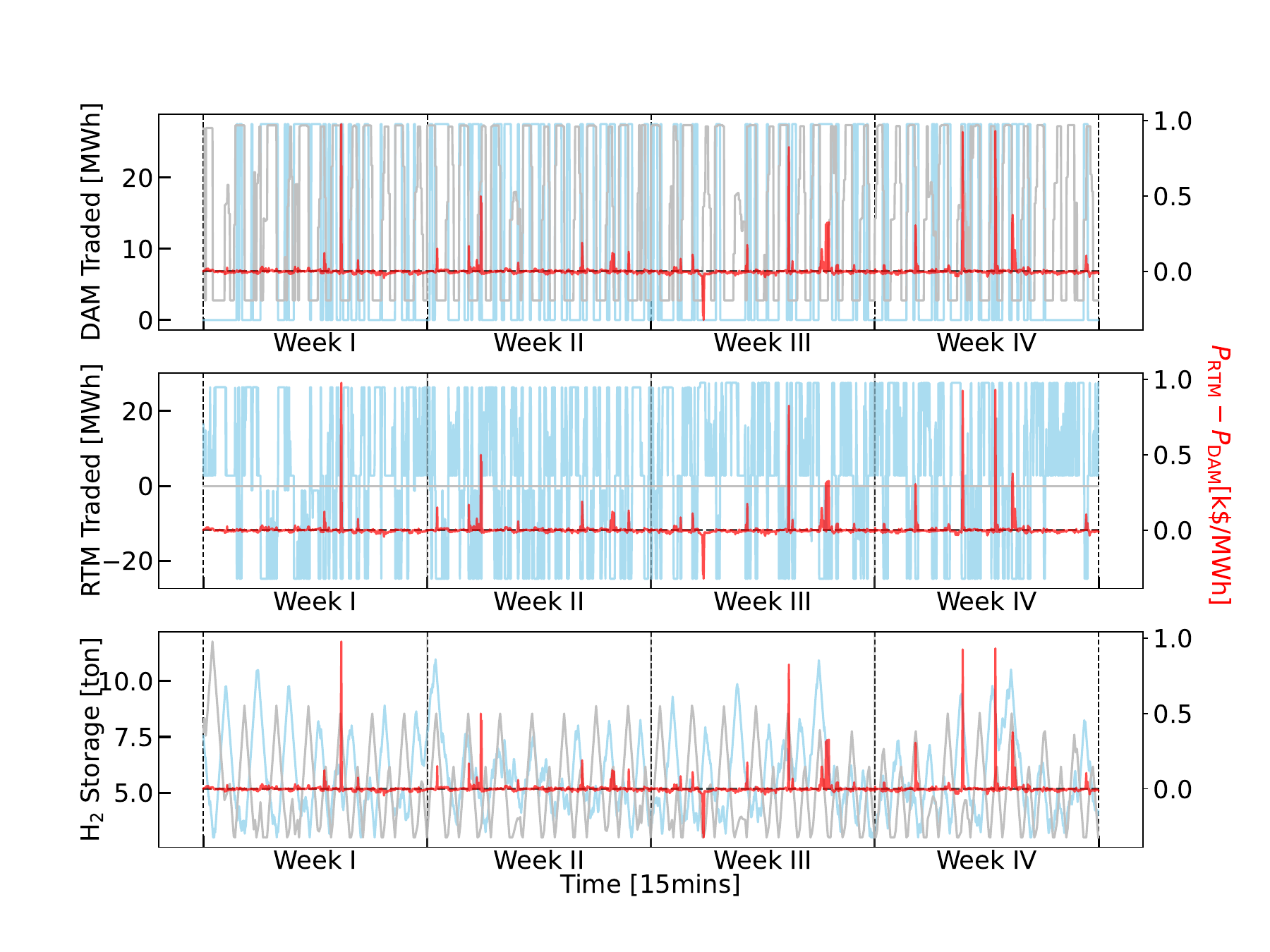}
        \caption{DAM traded power, RTM traded power, and $\text{H}_{\text{2}}$ storage levels for the HF-MS (light blue line) and LF-MS (gray line) strategies. Data correspond to four weeks in Houston, January 2022. The red line indicates the price differential (RTM Price minus DAM Price).}
    \label{fig:market_simple_ms_ss_cali}
\end{figure}

To further validate the proposed approach, we present additional experimental results for the California market. Similar to the Houston case, we compare the performance of the HF-MS and HF-SS strategies over four representative weeks, as shown in \Cref{fig:market_simple_ms_ss_cali}. The figure illustrates the traded DAM power (HF-MS in blue, HF-SS in gray), traded RTM power (HF-MS in blue, HF-SS in gray), hydrogen storage levels (HF-MS in blue, HF-SS in gray), and the RTM-DAM price difference (red line) from top to bottom. Several key observations can be made. First, in the DAM (top subplot), HF-MS strategically commits to significant power purchases, particularly in anticipation of high RTM prices (red line), indicating a proactive approach to securing lower-cost electricity. Second, HF-MS leverages these DAM purchases to capitalize on favorable RTM conditions (middle subplot). During periods of high RTM prices, HF-MS sells power back to the RTM (negative RTM traded values), profiting from the DAM-RTM price arbitrage. For instance, clear instances of RTM power sales during price spikes are evident in Week I and Week IV. Finally, the hydrogen storage dynamics (bottom subplot) highlight the operational flexibility of the HF-MS strategy. The frequent and substantial fluctuations in $\text{H}_{\text{2}}$ storage levels demonstrate the capacity of HF-MS to dynamically adjust hydrogen production in response to real-time market conditions and energy arbitrage opportunities, enabling more intensive utilization of the electrolysis plant.

\begin{figure}[!h]
 \centering
    \includegraphics[width=0.95\textwidth]{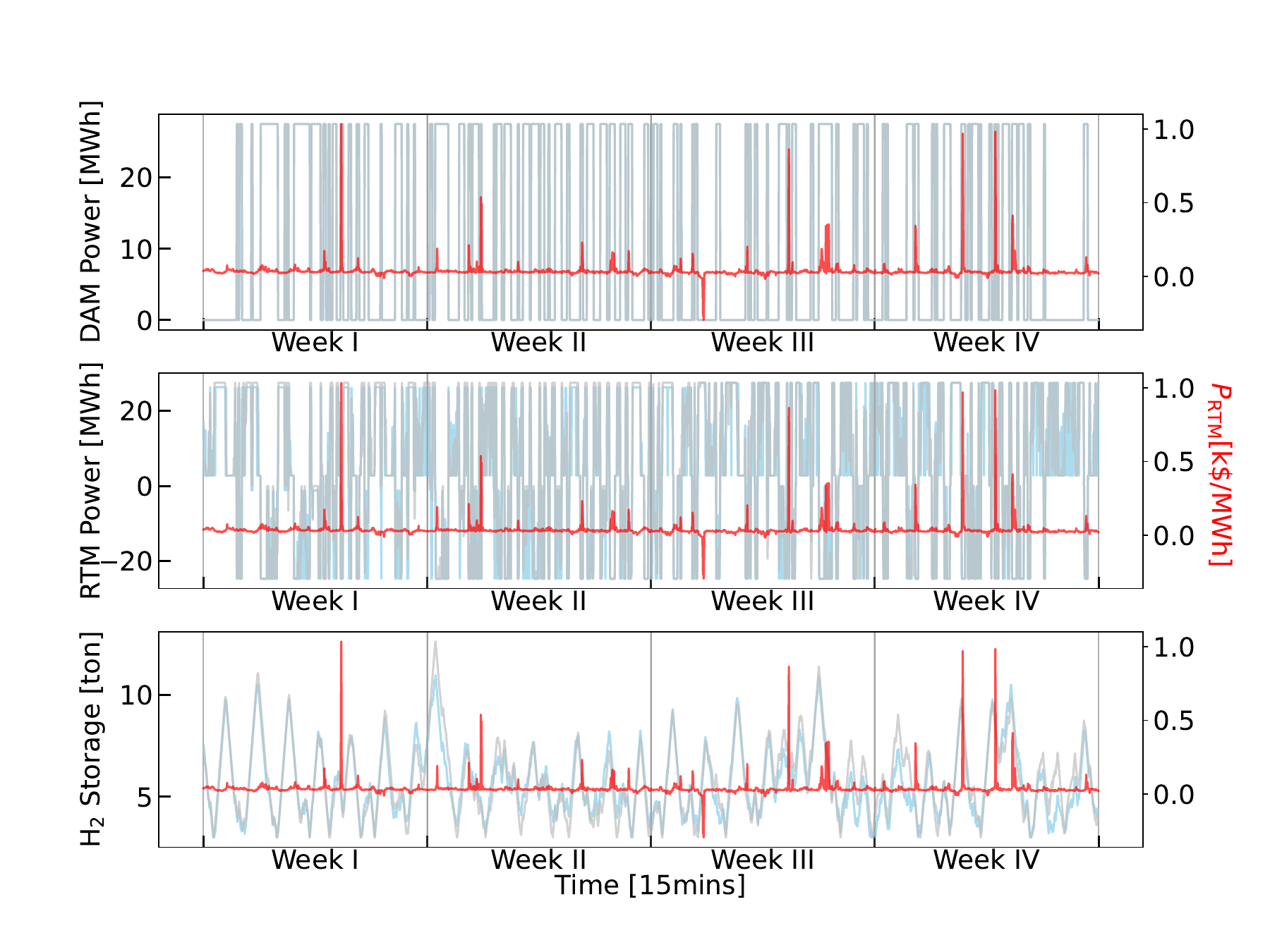}
        \caption{DAM traded power, RTM traded power, and $\text{H}_{\text{2}}$ storage levels for the HF-MS (light blue line) and LF-MS (gray line) strategies. Data correspond to four weeks in Houston, January 2022. The red line indicates the price differential (RTM Price minus DAM Price).}
    \label{fig:market_simple_hf_lf_cali}
\end{figure}

We further examine the results of HF-MS and LF-MS, as illustrated in \Cref{fig:market_simple_hf_lf_cali} for the California market. The key operational metrics for HF-MS are still represented by the light blue lines in \Cref{fig:market_simple_hf_lf_cali}. Notably, similar to the operational strategy observed in Houston, both approaches may utilize analogous foundational logics when engaging with the Day-Ahead Market (DAM) and Real-Time Market (RTM), primarily focusing on exploiting price differentials for arbitrage opportunities. However, the specific interaction patterns observed within the RTM—including the magnitudes of transactions and their consequent systemic impacts—reveal significant distinctions between the two strategies.

A key observation is the notably heightened engagement with the Real-Time Market (RTM) demonstrated by the LF-MS strategy. The middle subplot of \Cref{fig:market_simple_hf_lf_cali}, depicting RTM power transactions, illustrates that when LF-MS (represented by the gray line) acquires power from the RTM (indicated by negative RTM power values), the peak magnitudes and volumes of these transactions tend to be significantly larger. For instance, during the decline in RTM prices from Week II and Week IV (as represented by the red line), the LF-MS strategy not only actively participates in purchasing power but also exhibits notably pronounced sales volumes (depth of the light gray line extending below zero), suggesting a tendency for larger-volume reverse power flows when market conditions are favorable.

Consequently, this more assertive interaction with the RTM under the LF-MS strategy correlates with considerably greater fluctuations in hydrogen storage levels (shown in the bottom subplot, gray line). The hydrogen (\( \text{H}_{\text{2}} \)) storage undergoes cycles of significant accumulation, indicating that LF-MS's increased purchasing activity in the RTM—often aimed at capitalizing on advantageous pricing for enhanced hydrogen production—directly translates to more frequent and larger-amplitude variations in stored hydrogen. These intensified procurement efforts lead to accelerated cycles of hydrogen generation and storage, reflecting a more dynamic and potentially more rigorous utilization of hydrogen storage capacity.


A salient observation is the comparatively intensified RTM engagement exhibited by the LF-MS strategy. The middle subplot of \Cref{fig:market_simple_hf_lf_cali}, which details RTM power transactions, clearly illustrates that when LF-MS (gray line) purchase power from the RTM (indicated by negative RTM Power values), the peak magnitudes and volumes of these sales are often more substantial. For instance, during the RTM price drop in Week I to Week IV (red line indicating RTM price), the LF-MS strategy not only participates in purchasing power but the extent of these sales (depth of the light gray line beyond zero) appears particularly pronounced, suggesting a propensity for larger-volume reverse power flows when profitable. Consequentially, this more aggressive RTM interaction under the LF-MS strategy correlates with markedly more pronounced fluctuations in hydrogen storage levels (bottom subplot, gray line). The $\text{H}_{\text{2}}$ storage undergoes cycles of more substantial accumulation. This behavior indicates that LF-MS's more pronounced power purchasing activity in the RTM—often undertaken to capitalize on favorable pricing for intensified hydrogen production—directly translates to more frequent and larger-amplitude variations in stored hydrogen. Such intensified procurement efforts, leading to accelerated hydrogen generation and storage cycles, reflect a more dynamic, and potentially more strenuous, utilization of the hydrogen storage capacity.

\begin{figure}[!h]
 \centering
    \subfigure[Cumulative cost of California under various operation strategies.]{
    \includegraphics[width=1.00\textwidth]{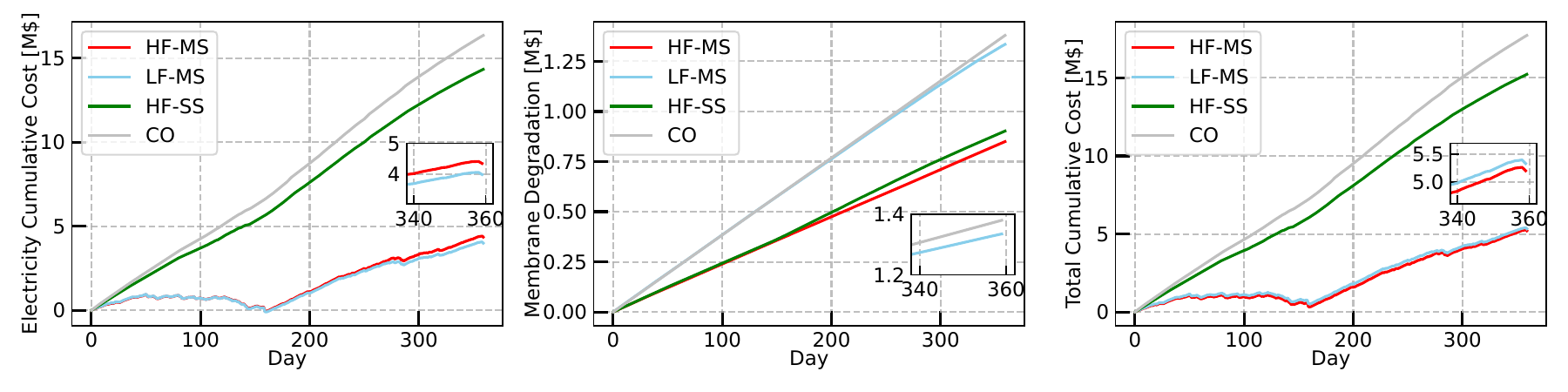}\label{subfig:cumCostbreakUpCali1}
} 
    \subfigure[Cumulative Cost Attribution for Different Operation Condition of California.]{
    \includegraphics[width=1.00\textwidth]{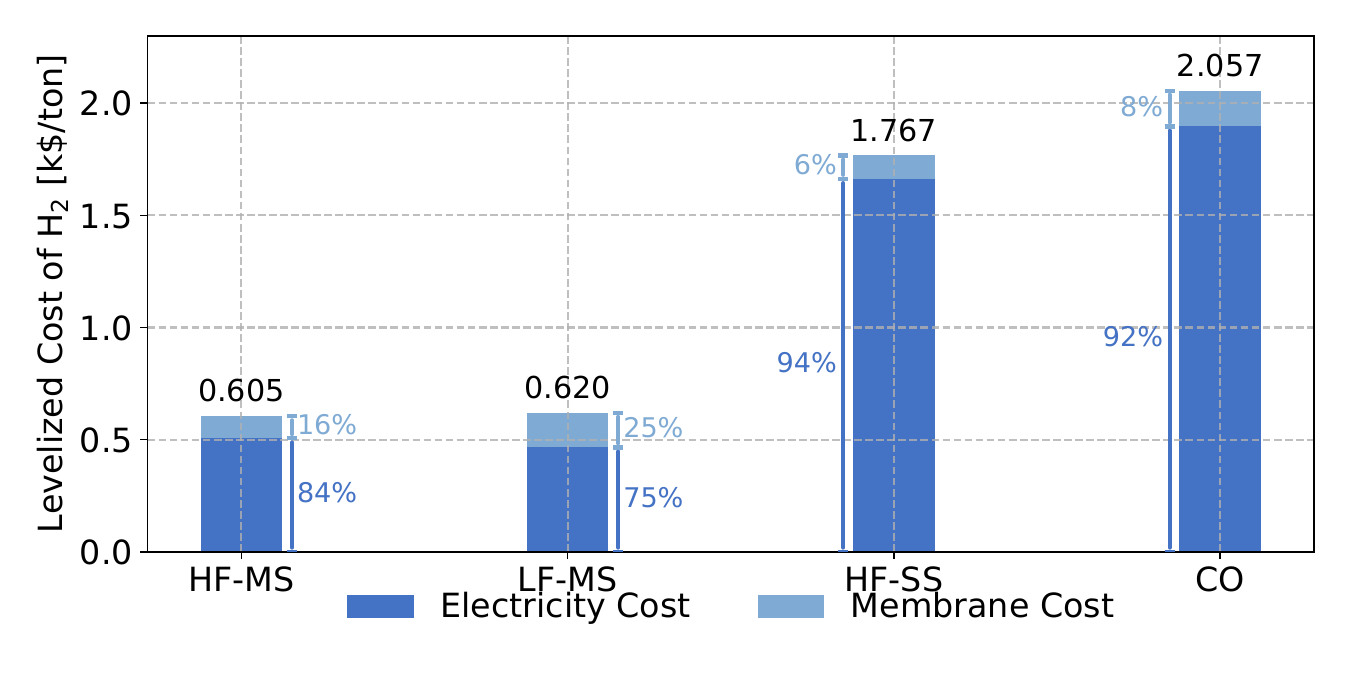}\label{subfig:cumCostbreakUpCali2}
} 
    \caption{Cost comparison results for California.}
    \label{fig:market_cost_compareCali}
\end{figure}

In addition, we compare the cumulative costs and levelized costs associated with different operational strategies for California, as illustrated in Figures \ref{fig:market_cost_compareCali}(a) and~\ref{fig:market_cost_compareCali}(b). We can draw parallels to our findings in Houston. Specifically, when operating within the multi-scale electricity market, the total cumulative cost is significantly lower than that of operational strategies that do not incorporate multi-scale considerations, as evidenced by the comparison of the red and light blue lines with the gray and green lines. It is noteworthy that the electricity price for HF-MS may be slightly higher than that of LF-MS; however, the membrane cost for HF-MS is lower than that for LF-MS. This phenomenon indicates that HF-MS tends to leverage membrane costs to offset electricity costs. Furthermore, in contrast to the Houston scenario, the California operations do not yield substantial profits, as the cumulative costs never fall below zero. This outcome can be attributed to the smaller price differentials between the RTM and DAM compared to Houston, leading to a limited profit margin. From Figure~\ref{fig:market_cost_compareCali}(b), we observe that the levelized cost of hydrogen for the multi-scale operation strategy is significantly lower than that of the single-scale operation strategies. This finding further underscores the necessity of considering trading between the RTM and DAM.

\begin{figure}[!h]
 \centering
    \subfigure[Operating Current Density and Temperature Vary Different Strategies in Houston.]{
    \includegraphics[width=0.8\textwidth]{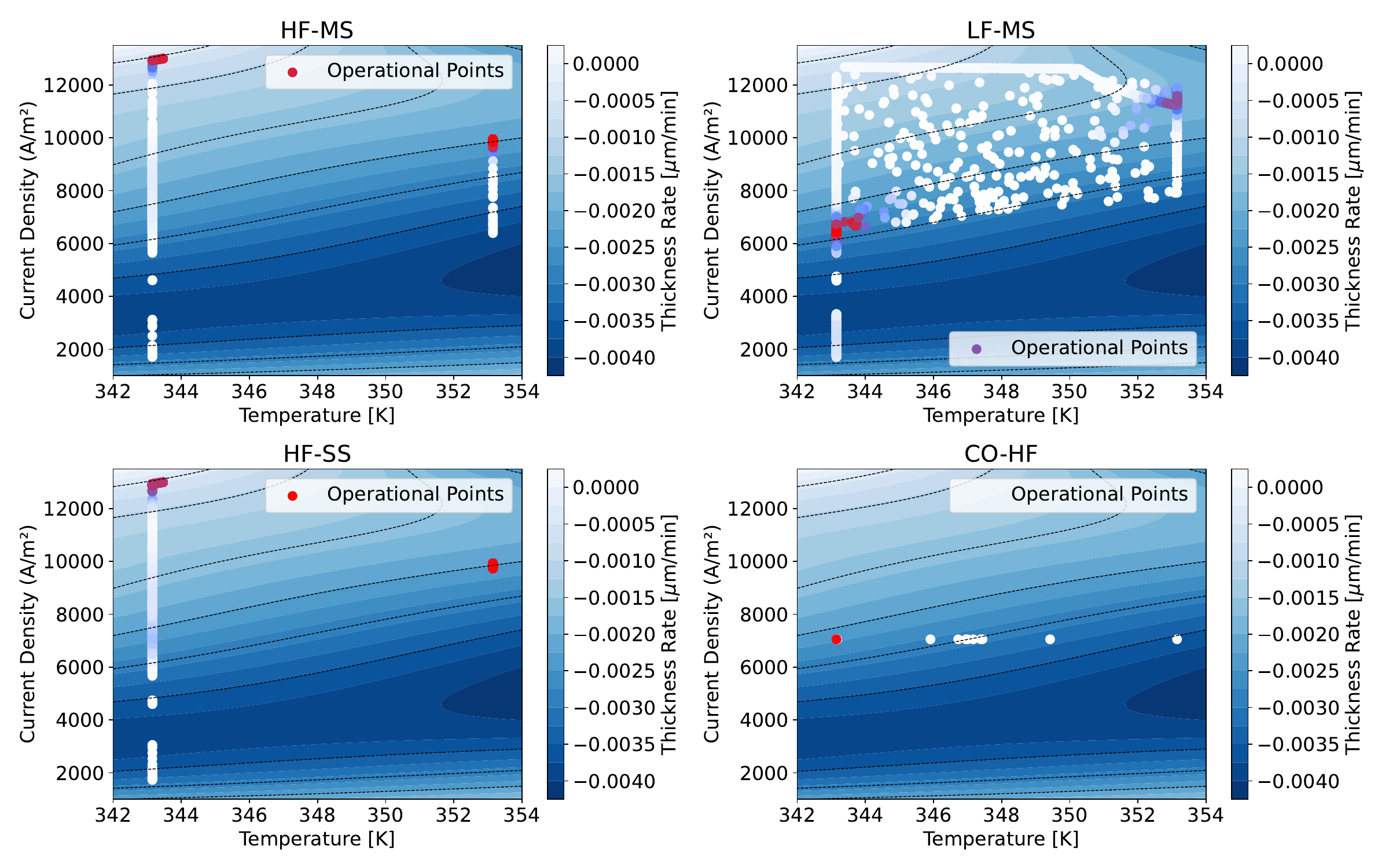} \label{fig:market_cost_breakdown1Cali}
    } 
    \subfigure[The KDE of current density at 343.15 K, the black line is membrane degradation rate, and the red line is the KDE estimation of operating current density.]{
    \includegraphics[width=1.0\textwidth]{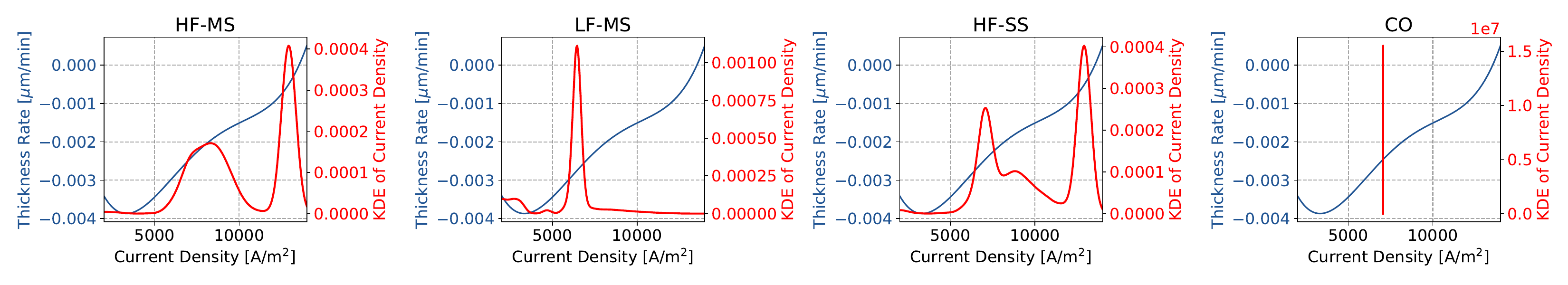} \label{fig:market_cost_breakdown2Cali}
    } 
    \subfigure[The KDE of current density at 353.15 K, the black line is membrane degradation rate, and the red line is the KDE estimation of operating current density.]{
    \includegraphics[width=1.0\textwidth]{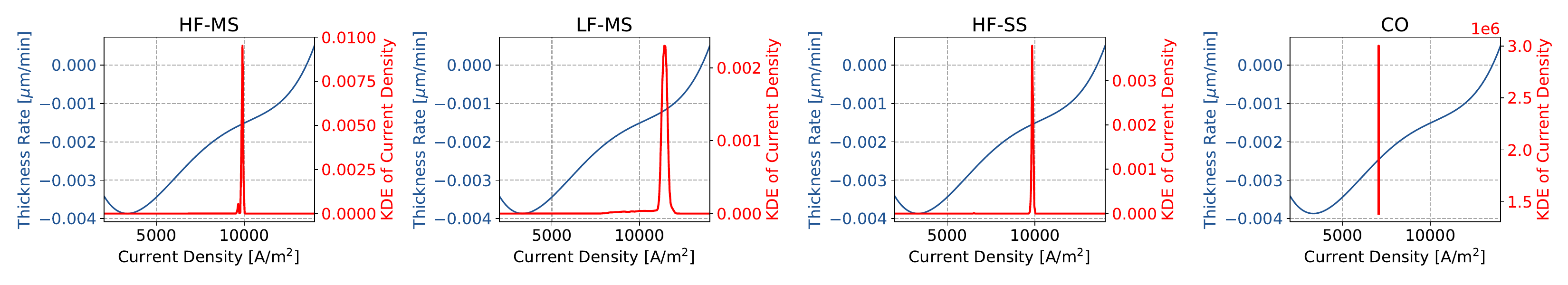} \label{fig:market_cost_breakdown3Cali}
    } 
    \caption{Membrane degradation rate and operating temperature in Califronia city.}
    \label{fig:market_cost_breakdownCali}
\end{figure}

Finally, we examined the operational factors contributing to the observed cost disparities. Our objective was to elucidate how the HF-MS and HF-SS strategies effectively mitigate membrane degradation. To achieve this, we analyzed the operating conditions of each strategy in relation to the membrane degradation rate, as illustrated in Figure~\ref{fig:market_cost_breakdownCali}(a).

From Figure~\ref{fig:market_cost_breakdownCali}(a), it is clear that both HF-MS and HF-SS consistently operate within regions characterized by lower temperatures and moderate current densities, thereby avoiding conditions that accelerate membrane degradation. In contrast, LF-MS and CO exhibit broader, more dispersed operational profiles, frequently venturing into high-risk zones for membrane failure.

Given that operations are concentrated at 343.15 K and 353.15 K, we further investigated the distribution of current densities and corresponding membrane degradation rates at these temperatures, as depicted in Figures~\ref{fig:market_cost_breakdownCali}(b) and \ref{fig:market_cost_breakdownCali}(c). The KDE results at 343.15 K (\Cref{fig:market_cost_breakdownCali}(b)) confirms that the HF-MS and HF-SS strategies cluster around current densities associated with minimal membrane degradation, resulting in a distinct advantage over LF-MS, with degradation rate differences on the order of 0.001 $\mu\text{m/min}$. These findings illustrate that high-fidelity operational control—especially when combined with multi-scale market participation—strategically minimizes exposure to conditions that exacerbate membrane failure, directly contributing to the observed reductions in cumulative membrane costs.

At 353.15 K, although this advantage is less pronounced, LF-MS exhibits slightly lower degradation rates compared to the HF strategies, with the difference being marginal (less than 0.0005 $\mu\text{m/min}$). This suggests that while higher operating temperatures modestly increase degradation rates for HF strategies, disciplined management of current density ensures that the overall impact remains small and comparable to LF strategies. Thus, the operational benefits of HF strategies are largely maintained even at elevated temperatures.

\end{appendices}

\ClearShipoutPicture
\end{document}